\newcommand{\cras}{0}
\newcommand{\faussepreuve}{0}

\newcommand{\definitE}{%
l'ensemble $\mathcal{E}$ des courbes à courbure algébrique positive, dont les extrémités et les tangentes en leurs extrémités sont données%
}
\newcommand{\definitEe}{%
the set $\mathcal{E}$ of curves with positive algebraic curvature, whose extremities and tangents in their extremities are given%
}

\newcommand{\resumefrancais}{%
On considère \definitE. \`A chacune des courbes de $\mathcal{E}$, on associe le minimum du rayon de courbure algébrique. 
\ifcase \cras
Nous montrons tout d'abord qu'il existe une courbe de $\mathcal{E}$ qui maximise ce minimum. 
Numériquement, on constate ensuite que cette courbe 
\or
\input{vrai_probleme_pointJ} 
\fi
est égale  à l'unique courbe de $\mathcal{E}$, formée d'un arc de cercle et d'un segment de droite, éventuellement réduit à un point. 
Cette courbe correspond aussi à un cas particulier des courbes de Dubins et sera 
utilisée pour améliorer la conception d'une pièce intervenant  dans un brevet.
}
\newcommand{\resumeanglais}{%
We consider \definitEe. 
For each of the curves of $\mathcal{E}$, we define the minimum of the radius of curvature. 
\ifcase \cras
We first prove that there exists a curve of $\mathcal{E}$ which maximizes this minimum.
Numerically, we observe then that this curve 
\or
\input{vrai_probleme_pointK} 
\fi
is equal  to the unique curve of $\mathcal{E}$ composed  of an  arc of circle 
and a line segment, 
where appropriate reduced to a point. This curve corresponds also to a particular case of Dubins's curve and will be used to improve the conception of a piece of a patent.
}

\newcommand{\adresse}{%
Laboratoire Inter-universitaire de Biologie de la Motricité\\
   POLYTECH\\
   Université Claude Bernard - Lyon 1\\
   15 Boulevard André LATARJET\\
   69622 Villeurbanne Cedex\\
France
}

\newcommand{\nom}{Jérôme Bastien}

\newcommand{\ftitrecras}{Existence et unicité d'une courbe à courbure positive maximisant le minimum du rayon de courbure}
\newcommand{\atitrecras}{Existence and uniqueness of a curve with positive curvature maximizing the minimum radius of curvature}
     \newcommand{\ftitre}{EXISTENCE D'UNE COURBE A COURBURE POSITIVE MAXIMISANT LE MINIMUM DU RAYON DE COURBURE -- "OBSERVATION NUM\'ERIQUE"}
\newcommand{\ftitrefacul}{EXISTENCE D'UNE COURBE MAXIMISANT LE MINIMUM DU RAYON DE COURBURE}

\newcommand{\aeL}{\text{p.p. sur  } ]0,L[}

\newcommand{\udsc}{\underline{ }}

\ifcase \cras

\documentclass[a4paper,11pt,reqno,final]{amsart}
\newcommand{\pathfilesty}{.}
\usepackage{\pathfilesty/articlefrancais}
\graphicspath{{./dessins/}}

\bibliography{%
../../../divers/circuit_rail/piece6_optimale/description_piece6_optimale,%
../../../divers/circuit_rail/piece6_optimale/description_piece6_optimale_maths,%
../../../divers/rapport_activite/mabibliographiebis,%
../../../divers/rapport_activite/monenseignementbib%
}

\begin{document}
\selectlanguage{french}



\title[\ftitrefacul]{\ftitre}

\author{\nom}

\date{\today}

\address{\adresse}

\email{\href{mailto:jerome.bastien@univ-lyon1.fr}{\nolinkurl{jerome.bastien@univ-lyon1.fr}}}

\begin{abstract}
\resumefrancais
\end{abstract}

\begin{abstract}
\selectlanguage{french}
\resumefrancais
\vspace{0.5 cm}

\noindent
\textsc{Abstract}.
\selectlanguage{english}
\resumeanglais
\selectlanguage{french}
\end{abstract}


\maketitle


\or


%


%

\documentclass{elsart3-1}


\usepackage{epsfig}

\usepackage{amssymb}

\usepackage[francais,english]{babel}

\newtheorem{theorem}{Theorem}[section]
\newtheorem{lemma}[theorem]{Lemma}
\newtheorem{e-proposition}[theorem]{Proposition}
\newtheorem{corollary}[theorem]{Corollary}
\newtheorem{e-definition}[theorem]{Definition\rm}
\newtheorem{remark}{\it Remark\/}
\newtheorem{example}{\it Example\/}
\newtheorem{theoreme}{Th\'eor\`eme}[section]
\newtheorem{lemme}[theoreme]{Lemme}
\newtheorem{proposition}[theoreme]{Proposition}
\newtheorem{corollaire}[theoreme]{Corollaire}
\newtheorem{definition}[theoreme]{D\'efinition\rm}
\newtheorem{remarque}{\it Remarque}
\newtheorem{exemple}{\it Exemple\/}
\renewcommand{\theequation}{\arabic{equation}}
\setcounter{equation}{0}

\usepackage[T1]{fontenc}
\usepackage[latin1]{inputenc}
\usepackage{url}
\usepackage[centerlast,small]{subfigure}
\usepackage{amsmath}
\usepackage{amsfonts}
\usepackage{psfrag}
\newenvironment{preuve}[1][]{\textsc{Preuve#1}}{\hfill \null \hfill  \qed}
\newenvironment{preuveidee}[1][]{\textsc{Idées de la preuve#1}}{\hfill \null \hfill \qed}
\newcommand{\Er}{\mathbb{R}}
\newcommand{\Erp}{\mathbb{R_{\, +}}}
\newcommand{\Erm}{\mathbb{R_{\, -}}}
\newcommand{\Erps}{\mathbb{R_{\, +} ^{\, *}}}
\newcommand{\Ce}{\mathbb{C}}
\newcommand{\En}{\mathbb{N}}
\newcommand{\norm}[2][{}]{\lVert#2\rVert_{{#1}}}
\newcommand{\vnorm}[2][{}]{\left\Vert#2\right\Vert_{#1}}
\newcommand{\prodsca}[2]{\left\langle #1,#2\right\rangle}  

\def\og{\leavevmode\raise.3ex\hbox{$\scriptscriptstyle\langle\!\langle$~}}
\def\fg{\leavevmode\raise.3ex\hbox{~$\!\scriptscriptstyle\,\rangle\!\rangle$}}

\journal{the Acad\'emie des sciences}
\begin{document}
\centerline{Géométrie/ Geometry}
\centerline{Problèmes mathématiques de la mécanique/Mathematical Problems in Mechanics}
\begin{frontmatter}

\selectlanguage{francais}
\title{\ftitrecras}

\author{\nom}
\address{\adresse}

\ead{jerome.bastien@univ-lyon1.fr}

\medskip
\selectlanguage{francais}

\begin{abstract}
\selectlanguage{francais}
\resumefrancais
\vskip 0.5\baselineskip

\selectlanguage{english}
\noindent{\bf Abstract}
\vskip 0.5\baselineskip
\noindent
{\bf  \atitrecras}
\resumeanglais
\end{abstract}
\end{frontmatter}

\graphicspath{{./dessins/}}


\fi


\selectlanguage{english}
\section*{Abridged English version}

In the framework of the patent \cite{brevetJB}, we had to define a curve 
whose  extremities $A$  et $B$  and the tangents at its points $A$ and $B$ are given,
both these  tangents  not being parallel.
The chosen curve is a parabola, or equivalently, a Bézier curve
(see also  
\cite{enumeration_circuit_JB_2016,%
bastienbrevetforum,%
bastienbrevetmmi,%
bastienbrevetapmep}). 
The disadvantage of this curve is that it has too small radii of curvature and we attempted  to find a  less incurved curve.
For this, we define \definitEe.
For each of curves of $\mathcal{E}$, we define the minimum of the radius of curvature.

This problem is very close to the problem of Dubins's curves \cite{MR0089457,MR0132460}, but it is not equivalent.
Dubins  also looks for curves whose extremities and tangents at the extremities are given. Under the assumption that the radius of curvature is everywhere
on the curve greater than a given radius of curvature $R>0$, he determines the curve which minimizes the length.
He proves that these curves, entirely defined by $R$, called geodesic and denoted
\begin{equation}
\label{rapeldubinseq01}
\mathcal{G}(R),
\end{equation}
are necessarily the union of three arcs of circle of radius $R$,
or the union of two arcs of circle of radius $R$ and of line segment or  a subset of these curves.
In our case, 
we impose the positive algebraic curvature and 
we do not consider \textit{a priori} this radius 
\ifcase \cras
of curvature $R$ and by adapting the proof of Dubins, we prove that there exists a curve of $\mathcal{E}$, which maximizes the minimum of algebraic radius of curvature.
\or
of curvature $R$.
\fi
We 
\ifcase \cras
secondly 
\fi
prove that there exists a unique curve of $\mathcal{E}$ composed of an arc of circle of radius $R_a$ and of a line segment, denoted $\mathcal{J}$.
This case corresponds to the limit case of figure \ref{coducl1}.
This radius $R_a$ depends only on the points $A$ and $B$ and on the tangents to these points and 
 is the greatest value of $R$, for which the Dubins's curve 
$\mathcal{G}(R)$ belongs to  $\mathcal{E}$.
\ifcase \cras
\or
\input{vrai_probleme_pointE}
\fi

\ifcase \cras
We consider then the set of the curves of $\mathcal{E}$ composed of a finite number $p\in \En^*$ of arcs of circles, 
each of them being of radius $R_k\in \Erps\cup \{+\infty\}$.
This set possesses also a curve which maximizes the minimum of radius of curvature, \textit{i.e.}, which maximizes the minimum of the $R_k$.
We then choose any $p \in \En^*$ and we assume that each arc of circle is defined by an angle $\theta_0$ constant
(the angle $\widehat{A_{k}\Omega{k}A_{k+1}}$, where $A_k$ and $A_{k+1}$ are the two extremities of the $k$-th arc of circle and $\Omega_k$ its center).
We observe that the optimal curve possesses two parts: on the first one, the arcs of circles have large radii of curvature, and on the second one,
the arcs of circles have radii which seem be close to the value of $R_a$  defined  above.
So, numerically, as $p$ tends to infinity, the  curve obtained seems to be close to  the unique curve 
  $\mathcal{J}$  of  $\mathcal{E}$. 

The parabola used in patent \cite{brevetJB} has been replaced by this curve and this allows us to obtain a new minimum radius of curvature
equal to $(\sqrt{2}-1)/2\approx 0.207$     instead of $1/25\,\sqrt {5} \approx 0.089$.

The important lack of this present 
\ifcase \cras
Article,
\or
Note, 
\fi
that we have to improve on, is to prove that which we observed numerically, \textit{i. e.} that the unique curve of  $\mathcal{E}$, composed of a line segment and of an arc of circle
is one of the curves which maximizes the minimum of radius of curvature. 
It would be interesting to prove the uniqueness of these curves, if it is true, in the continuous and the discrete cases.

\fi

\ifcase \cras
\selectlanguage{french}
\or
\selectlanguage{francais}
\fi


\ifcase \cras

\section*{Version abrégée en français}

Dans le cadre du brevet  \cite{brevetJB}, il a été nécessaire de construire une courbe
passant par deux points $A$ et $B$ du plan, dont les directions des tangentes en $A$  et $B$ 
sont données en étant non parallèles. La courbe choisie est une parabole, ou de façon équivalente, une courbe de Bézier 
(voir aussi 
\cite{enumeration_circuit_JB_2016,%
bastienbrevetforum,%
bastienbrevetmmi,%
bastienbrevetapmep}). 
Cette courbe présentant l'inconvénient d'avoir des rayons de courbures trop petits, on a essayé de trouver une courbe moins incurvée. Pour cela, 
on définit \definitE\ et à chacune des courbes de $\mathcal{E}$, on associe le minimum du rayon de courbure.
\ifcase \cras
\or
On montre tout d'abord qu'il existe une courbe de $\mathcal{E}$ qui maximise ce minimum.
\fi

Ce problème ressemble fortement aux courbes de Dubins \cite{MR0089457,MR0132460}
sans lui être équivalent.
Dubins cherche des courbes  passant aussi par deux points $A$ et $B$ du plan, dont les directions des tangentes en $A$ en $B$ 
sont données. Sous l'hypothèse qu'en tout point, le rayon de courbure est supérieur à $R$ où $R>0$ est donné à l'avance, il cherche la courbe qui minimise la longueur.
Il montre que de telles courbes, entièrement définie par $R$ appelées  géodésiques et notées \eqref{rapeldubinseq01},
ne peuvent être que la réunion de trois arcs de cercles de rayon $R$
ou la réunion de deux  arcs de cercle de rayon $R$ et d'un segment, ou une sous partie de ces deux courbes.
Dans notre cas, 
nous imposons la courbure algébrique positive,
nous ne donnons pas ce rayon de courbure $R$ \textit{a priori}.
\ifcase \cras
En adaptant la démonstration de Dubins, on montre qu'il existe une courbe de $\mathcal{E}$, maximisant le minimum du rayon de courbure algébrique,
dont on impose qu'il reste positif.
\fi
Nous montrons 
\ifcase \cras
ensuite
\or
tout d'abord 
\fi
qu'il existe une unique courbe de $\mathcal{E}$ formée d'un arc de cercle de rayon $R_a$ et d'un segment de droite,
notée $\mathcal{J}$.
Ce cas   correspond au cas limite de la figure \ref{coducl1}.
Ce rayon $R_a$  dépend uniquement des points $A$ et $B$ et des tangentes données en ces points et  est la plus grande valeur possible de $R$, pour laquelle 
la courbe de Dubins
$\mathcal{G}(R)$ est dans $\mathcal{E}$.
\ifcase \cras
\or
\input{vrai_probleme_pointF}
\fi

\ifcase \cras
Ensuite, on considère l'ensemble des courbes de $\mathcal{E}$ formée un nombre de fini $p\in \En^*$ d'arcs de cercles, chacun de rayon $R_k\in \Erps\cup \{+\infty\}$.
Cet ensemble possède une courbe qui maximise le minimum du rayon de courbure, c'est-à-dire, qui maximise le minimum des $R_k$.
On choisit ensuite $p$ quelconque et on impose que chaque arc de cercle
soit défini par un angle $\theta_0$ constant (l'angle $\widehat{A_{k}\Omega{k}A_{k+1}}$, où $A_k$ et $A_{k+1}$ sont les deux extrémités du $k$-ième  l'arc de cercle et $\Omega_k$ son centre).
On constate que la courbe optimale obtenue possède deux parties : sur l'une d'elle, les arcs de cercles ont un grand rayon de courbure, sur l'autre,
les arcs de cercles ont un rayon qui semble se rapprocher de la valeur $R_a$ définie ci-dessus. Ainsi, numériquement, quand $p$ grandit la courbe obtenue semble se rapprocher de
l'unique courbe $\mathcal{J}$  de $\mathcal{E}$. 
\fi

\ifcase \cras
L'importante lacune de 
\ifcase \cras
cet article,
\or
cette note, 
\fi
qu'il conviendrait de combler  par la suite, est de montrer ce que l'on a observé numériquement, c'est-à-dire que 
l'unique courbe de $\mathcal{E}$, formée d'un segment de droite et d'un arc de cercle de rayon $R_a$ est bien l'une des courbes
qui maximise le minimum du rayon de courbure.
Par la suite, il serait intéressant de montrer les éventuelles unicités des courbes, dans le cas continu comme dans les différents cas discrets évoqués.
\fi

\fi

\section{Introduction}
\label{introduction}

Dans le cadre du brevet  \cite{brevetJB}, il a été nécessaire de construire six courbes 
de classe ${\mathcal{C}}^1$  chacune d'elles passant par un point $A$ et un point $B$ en étant tangente 
respectivement en $A$ et $B$ aux droites $(AO)$ et $(BO)$, où, pour chacune d'elle, $A$, $B$ et $O$ sont trois points donnés du plan.
Plus de détails pourront être trouvés dans 
\cite{enumeration_circuit_JB_2016,%
bastienbrevetforum,%
bastienbrevetmmi,%
bastienbrevetapmep}. Chacune de ces courbes doit relier un des sommets ou un des milieux de côté d'un carré de centre
$O$ et de côté $1$, le point $O$ est fixé, centre du carré et conventionnellement choisi comme repère,  et les points $A$ et $B$ sont l'un des sommets ou un des milieux de côté du carré.
Compte tenu des isométries laissant invariant le carré, seules six courbes ont dû être définies : 
deux segments de droites, deux arcs de cercles et deux portions de paraboles, comme représenté sur 
\cite[Figure 1]{enumeration_circuit_JB_2016}.
On définit les deux vecteurs $\vec \alpha$ et $\vec \beta$  et l'angle $\Omega$ de la façon suivante : 
\begin{subequations}
\label{eq10tottot}
\begin{align}
\label{eq10tot}
&\vec \alpha=\frac{1}{OA}\overrightarrow{AO},\quad 
\vec \beta=\frac{1}{OB}\overrightarrow{OB},
\intertext{et}
\label{eq10totnew}
&\Omega=
\left( \widehat{\vec \alpha,\vec \beta} \right).
\end{align}
\end{subequations}
De façon plus générale, on 
se donne trois du plan, $O$,  $A$ et $B$, deux à deux distincts, deux vecteurs unitaires $\vec \alpha$ et $\vec \beta$ et l'angle $\Omega$ définis par
\eqref{eq10tottot},
où $\Omega$ n'est pas un multiple de $\pi$. 
Quitte 
à changer de sens de parcours de la courbe, donc à intervertir $A$ et $B$ et à  multiplier $\vec \alpha$ et $\vec \beta$ par $-1$, 
 on peut supposer, sans perte de généralité, que $\Omega$ vérifie 
\begin{equation}
\label{eq20}
\Omega \in ]0,\pi[.
\end{equation}
On s'intéresse à une courbe au moins de classe ${\mathcal{C}}^1$ passant par $A$, tangente à $(OA)$ en $A$, tangente à 
$(OB)$ en $B$. 
On peut choisir une parabole, ce qui a été fait par exemple dans le cas du brevet sur 
les \cite[figures 1e) et 1f)]{enumeration_circuit_JB_2016}

Ces courbes ont servi à définir des rails \textit{Easyloop} \textregistered\ aptes à faire rouler un train miniature.
Lors de la fabrication des pièces,
la dernière forme, donnée dans 
\cite[La  figure 1f)]{enumeration_circuit_JB_2016}
n'a pas été retenue, puisque trop incurvée, c'est-à-dire, 
que le minimum du rayon de courbure  est trop petit.
Il est en effet nécessaire que le rayon de courbure ne soit pas trop petit pour deux raisons. 
Tout d'abord, les courbes ainsi définies correspondent aux lignes médianes des rails construits : 
les passages des roues et les bords du rails sont définis comme des courbes à distance constante de ces courbes et si le rayon de courbure est trop petit, des points stationnaires avec
des discontinuités de la tangente peuvent apparaître. En outre, les roues des véhicules qui empruntent ces rails doivent pouvoir tourner par rapport au châssis du véhicule et si le rayon de courbure
est trop petit, l'angle de braquage, c'est-à-dire, l'angle entre les essieux qui supportent les paires de roues et l'axe longitudinal du véhicule, est trop important.
Nous proposons donc de déterminer une courbe pour définir la pièce correspondant à 
\cite[figure 1f)]{enumeration_circuit_JB_2016}
qui soit optimale, au sens où le minimum du rayon de courbure est choisi le plus grand possible. 
Nous imposerons aussi à la courbe recherchée d'être à courbure positive. Sans cette hypothèse, le problème est mal posé, puisque
l'on peut construire une courbe formée de trois arcs de cercles, chacun de rayon $R$, avec $R$ arbitrairement grand.

Ce problème ressemble fortement à celui des courbes de Dubins \cite{MR0089457,MR0132460}
sans lui être équivalent.
Dubins cherche des courbes  passant aussi par deux points $A$ et $B$ du plan, dont les directions des tangentes en $A$ en $B$ 
sont données. Sous l'hypothèse qu'en tout point, le rayon de courbure est supérieur à $R$ où $R>0$ est donné à l'avance, il cherche la courbe qui minimise la longueur.
Il montre que de telles courbes, entièrement définie par $R$ appelées  géodésiques et notées 
\eqref{rapeldubinseq01}
ne peuvent être que la réunion de trois arcs de cercles de rayon $R$
ou la réunion de deux  arcs de cercle de rayon $R$ et d'un segment, ou une sous partie de ces deux courbes.
\newcommand{\lengthdiftypcodu}{6.5}
\begin{figure}
\psfrag{A}{$A$}
\psfrag{B}{$B$}
\psfrag{O}{$O$}
\psfrag{C}{}
\psfrag{D}{}
\psfrag{E}{}
\psfrag{F}{}
\psfrag{G}{}
\psfrag{a}{$\vec \alpha$}
\psfrag{b}{$\vec \beta$}
\centering
\subfigure[\label{codu0}Le cas $0$ : $0<R<R_a$\ifcase \cras \or /\textit{{The case $0$ : $0<R<R_a$}}\fi]
{\epsfig{file=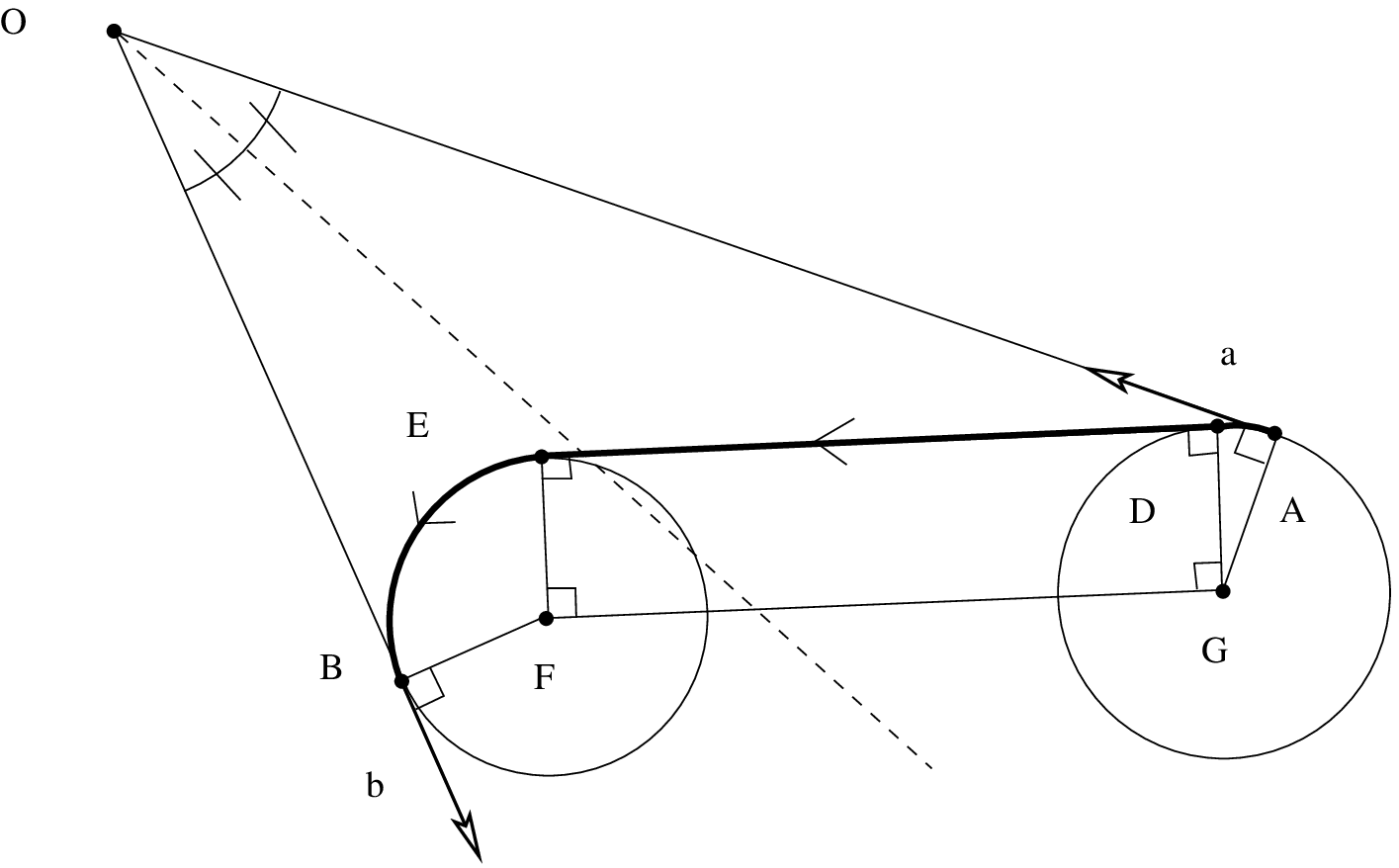, width=\lengthdiftypcodu cm}}
\quad
\subfigure[\label{coducl1}Le cas limite $0-1$ : $R=R_a$\ifcase \cras \or //\textit{{The limit case $0-1$ : $R=R_a$}}\fi]
{\epsfig{file=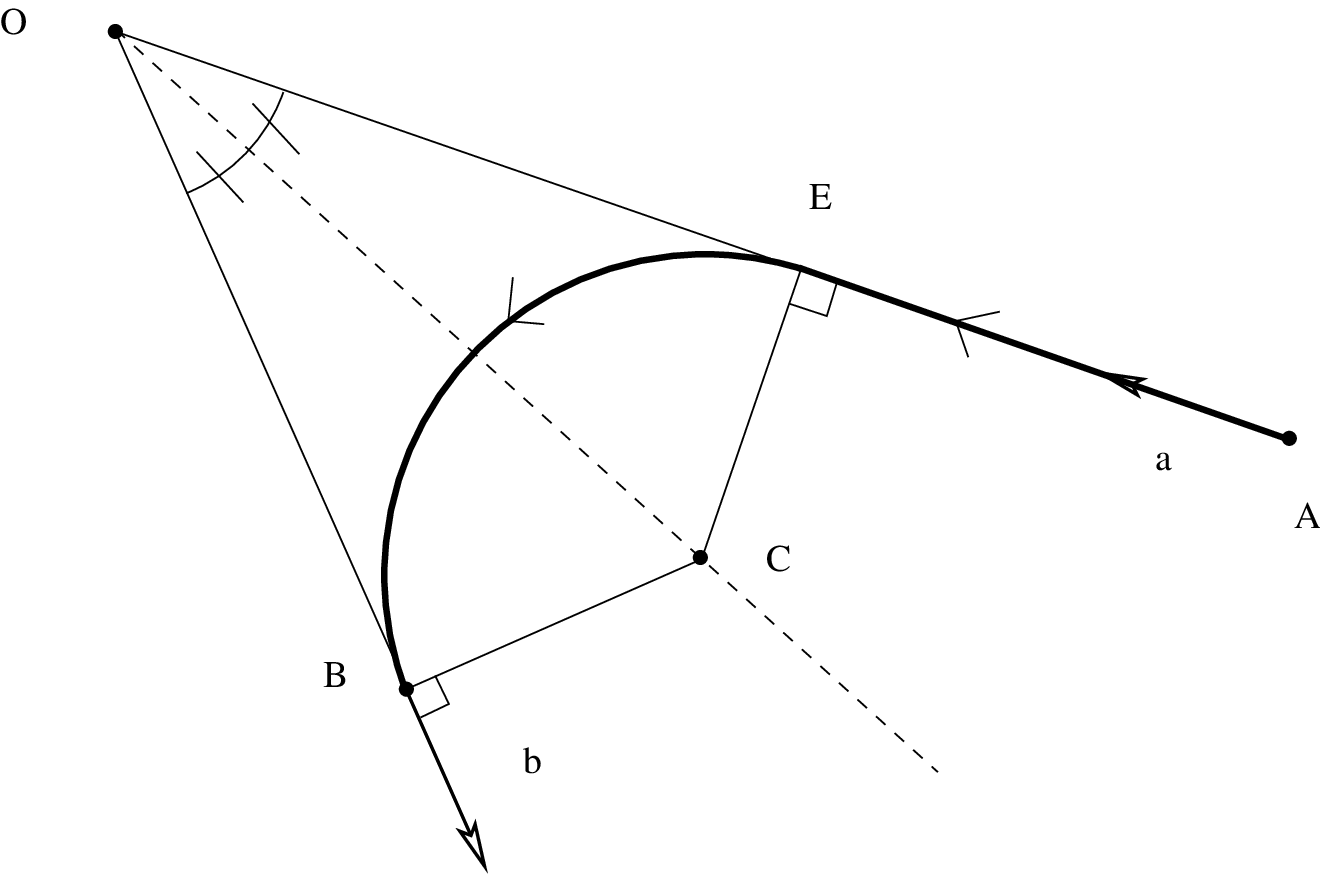, width=\lengthdiftypcodu  cm}}
\caption{\label{diftypcodu}Les différents types de courbes de Dubins définies par $R$ correspondant aux points $A$, $B$ et les vecteurs $\alpha$ et $\beta$ dans les cas $0<R\leq R_a$.\ifcase \cras \or/\textit{The different kinds of Dubins's curves 
defined by $R$ corresponding to points $A$ and $B$ and to vectors  $\alpha$ and $\beta$ in the cases  $0<R\leq R_a$.}\fi}
\end{figure}
Les courbes de Dubins correspondant à $0<R\leq R_a$ (où le rayon $R_a$ ne dépend que des points $O$, $A$ et $B$) sont représentées sur la figure \ref{diftypcodu}.
On considère \definitE. 
Dans notre cas, à la différence des travaux de Dubins, nous ne donnons pas ce rayon de courbure $R$ \textit{a priori} et nous imposons une courbure algébrique positive. 
\ifcase \cras
En adaptant la démonstration de Dubins, on montre qu'il existe une courbe de $\mathcal{E}$, maximisant le minimum du rayon de courbure.
\fi
Nous montrerons qu'il existe une unique courbe de $\mathcal{E}$ formée d'un arc de cercle de rayon $R_a$ et d'un segment de droite.
Ce cas   correspond au cas limite de la figure \ref{coducl1}.
Notre problème, qui ne me semble pas évoqué dans la littérature\footnote{Cette question 
a néanmoins été partiellement soulevée, mais visiblement non résolue sur 
\url{https://math.stackexchange.com/questions/1391778/connect-two-points-given-their-angles-with-a-maximum-radius}},
est donc distinct \textit{a priori} de celui de Dubins.
Ces travaux de Dubins ont été retrouvés plus tard autrement en utilisant le principe de maximum de Pontryagin par exemple dans 
\cite{Boissonnatetal1991,Boissonnatetal1994}. 
Très utilisées en robotique et en automatique, ces courbes de Dubins ont fait l'objet de nombreuses recherches.
Voir par exemple les deux thèses suivantes
\cite{lazardtel00442770,Jalel2016}.
De nombreuses variantes existent sur ces courbes de Dubins : 
si des point de rebroussement sont possibles (ce que l'on n'a pas ici, puisque le paramétrage est normal) dans le cas 
où le robot peut inverser sa vitesse \cite{MR1069892} ; 
les recherches prenant en compte les obstacles ont été initiées par Laumond dans 
\cite{Laumond1987};  des problèmes analogues avec des courbes constituées d'arcs de cercle et de segments de longueurs minimales imposées
sont présentées dans 
\cite{MR3158582}.
Notons enfin que dans \cite{Moser2009}, un problème proche de  notre problème est évoqué : il s'agit de trouver la
courbe, de longueur donnée (ou inférieure à une longueur donnée) 
qui maximise le minimum du rayon de courbure divisé par le rayon de courbure d'une courbe de référence, donnée à l'avance, 
comme frontière d'un convexe donné.

Dans la section \ref{enonce}, nous énonçons le problème. 
\ifcase \cras
\ifcase \faussepreuve
Dans la section \ref{existencetheonewjuste},  nous montrerons qu'il existe une courbe de $\mathcal{E}$, maximisant le minimum du rayon de courbure, c'est-à-dire 
minimisant  le maximum de la courbure.
\or
Dans la section \ref{existencetheo}, nous montrerons qu'il existe une courbe de $\mathcal{E}$, maximisant le minimum du rayon de courbure, c'est-à-dire 
minimisant  le maximum de la courbure, \textbf{résultat abandonné, car malheureusement faux ! Consulter cette section ! Rétabli !!!!}
\fi
\fi
En section \ref{uniquecercledroite}, nous construirons  l'unique courbe de $\mathcal{E}$ formée d'un segment de droite et d'un arc de cercle 
et nous montrerons ensuite  que cette courbe correspond à une courbe 
 optimale parmi les courbes de Dubins. 
\ifcase \cras
Enfin, en section \ref{discret}, 
nous présenterons le problème discret et donnerons une simulation numérique 
qui semble faire apparaître le fait que l'unique  courbe de $\mathcal{E}$ formée d'un segment de droite et d'un arc de cercle est bien l'une de celles
qui maximisent le minimum du rayon de courbure.
\or
Le \cite[théorème 3.2]{piece6_optimale_JB_2019_X4}  affirme  
qu'il  existe une courbe $X$ de $\mathcal{E}$ qui minimise le maximum de la courbure.
De plus, des simulations numériques présentées dans \cite[section 5]{piece6_optimale_JB_2019_X4} ont corroboré le fait que 
l'unique  courbe de $\mathcal{E}$ formée d'un segment de droite et d'un arc de cercle est bien l'une de celles
qui maximisent le minimum du rayon de courbure.
\input{vrai_probleme_pointA}
\fi

\section{\'Enoncé du problème}
\label{enonce}

Reprenons le formalisme du papier historique de Dubins  \cite{MR0089457}.
On se donne $O$, $A$ et $B$ trois points deux à deux distincts du plan et $\vec \alpha$ et $\vec \beta$ 
deux vecteurs, vérifiant 
 \eqref{eq10tottot} et  \eqref{eq20}.
Nous cherchons une courbe, paramétrée de façon normale par son abscisse curviligne : 
$s\in [0,L=L(X)]\mapsto X(s)\in \Er^2$, de classe ${\mathcal{C}}^1$.
Pour toute la suite, pour toute telle fonction $X$ de classe  ${\mathcal{C}}^1$, on notera 
\begin{equation*}
L=L(X).
\end{equation*}
On note $\vnorm{.}$ la norme Euclidienne de $\Er^2$. 
On suppose que l'on a 
\begin{subequations}
\label{eq100tot}
\begin{align}
\label{eq100a}
&\forall s\in [0,L],\quad \vnorm{X'(s)}=1,\\
\label{eq100b}
&X(0)=A,\\
\label{eq100c}
&X(L)=B,\\
\label{eq100d}
&X'(0)=\vec \alpha,\\
\label{eq100e}
&X'(L)=\vec \beta.
\end{align}
\end{subequations}
On supposera que 
\begin{equation}
\label{eq110}
X\in W^{2,\infty}(0,L;\Er^2),
\end{equation}
ce qui permet de définir  la valeur absolue de la courbure  $c$  par 
\begin{equation}
\label{eq115old}
|c(s)|= \vnorm{X''(s)}, \quad \aeL.
\end{equation}
La fonction $X$ est dans $ {\mathcal{C}}^1([0,L];\Er^2)$, on a $X('s)\not =0$ 
et on considère une détermination continue de l'angle 
$\phi$ défini par 
\begin{equation}
\label{eq120}
\forall s\in [0,L],\quad 
 \phi(s)=
\left( \widehat{\vec \alpha,X'(s)} \right).
\end{equation}
La fonction $\phi$ est donc dans 
$W^{1,\infty}(0,L)\subset  {\mathcal{C}}^0([0,L])$
et on a 
\begin{equation}
\label{eq130}
\frac{d \phi}{ds}=c,
 \quad \aeL.
\end{equation}
où ici $c$ désigne la courbure algébrique. 
On impose alors 
\begin{equation}
\label{eq140}
\text{$c$ est presque partout positive.}
\end{equation}
Ainsi, \eqref{eq140}
 est équivalent à 
\begin{equation}
\label{eq150}
\text{$\phi$ est croissant.}
\end{equation}
Dans ce cas, on peut réécrire \eqref{eq115old} sous la forme 
\begin{equation}
\label{eq115}
c(s)= \vnorm{X''(s)}, \quad \aeL.
\end{equation}
Notons que \eqref{eq150} implique 
la condition suivante :
\begin{equation}
\label{eq160}
\forall s\in [0,L],\quad 
 \phi(s)\in [0,\Omega].
\end{equation}
Dire 
que le minimum du rayon de courbure est le plus grand possible revient à dire
que le maximum de la valeur absolue de la courbure est  le plus petit possible, soit encore, selon 
\eqref{eq140},  que le maximum de la courbure, définit comme $\vnorm[L^\infty(0,L)]{c}$  est  le  plus petit possible.

\begin{definition}
\label{defiE}
Pour $O$, $A$ et $B$ trois points du plan deux à deux distincts et $\vec \alpha$ et $\vec \beta$ 
deux vecteurs vérifiant 
 \eqref{eq10tottot} et \eqref{eq20}, on définit l'ensemble $\mathcal{E}$ des courbes $X$ du plan 
vérifiant 
\eqref{eq100tot},
\eqref{eq110},
\eqref{eq140} (ou \eqref{eq150}).
\end{definition}

D'après \eqref{eq115}, le problème consistera finalement à déterminer une courbe $X$ de $\mathcal{E}$ minimisant 
$\vnorm[L^{\infty}(0,L)]{\vnorm{X''}}$, c'est-à-dire le sup essentiel de la fonction de $[0,L]$ dans $\Er$ : 
$s\mapsto \vnorm{X''}$ : 
\begin{equation}
\label{eq200}
\vnorm[L^{\infty}(0,L)]{\vnorm{X''}}
=\min_{Z\in \mathcal{E}} \vnorm[L^{\infty}(0,L)]{\vnorm{Z''}}.
\end{equation}

Notons enfin que $X$ est de classe ${\mathcal{C}}^1$ et non nécessairement ${\mathcal{C}}^2$.
D'autres travaux  utilisent  
les courbes à courbures continues, utilisant par exemples les clothoïdes comme courbes de raccordement, 
qui permettent une croissance continue de la courbure \cite{MR1717119,scheuertel00001746}. 
En effet, 
dans le cadre de la robotique ou du transport,
il n'est pas possible d'avoir une discontinuité de la courbure, qui implique une discontinuité des accélérations normales, et 
donc des chocs, ce qui use le matériel prématurément ou gêne le voyageur ; un robot ou un véhicule ne peut pas non plus changer instantanément d'angle de braquage.
Au contraire ici,  la discontinuité de la courbure  ne nous gêne pas pour plusieurs raisons. 
Dans le domaine du jouet, 
les masses et les vitesses des véhicules sont très faibles, donc les chocs dus aux discontinuité de l'accélération normale sont négligeables. 
De plus, la notion de confort du voyageur n'a pas de sens.
Les roues des véhicules peuvent subir une discontinuité de l'angle de braquage 
parce qu'elles présentent un léger jeu par rapport au châssis. 
Enfin, la courbe construite dans le cas du brevet 
 \cite{brevetJB,enumeration_circuit_JB_2016} est de classe ${\mathcal{C}}^1$, ${\mathcal{C}}^2$ par morceaux,
mais non ${\mathcal{C}}^2$, puisque formée de portions de segments, de cercles et de paraboles. 
Il n'est donc pas nécessaire de restreindre notre étude aux courbes ${\mathcal{C}}^2$.

\ifcase \cras
\begin{remark}
\or
\begin{remarque}
\fi
Quitte à parcourir, le cas échéant, la courbe dans l'autre sens, on peut 
remplacer respectivement
\eqref{eq20} par "$\Omega \in ]-\pi,\pi[\setminus \{0\}$", 
\eqref{eq140} par "$c$ est presque partout de signe constant", 
\eqref{eq150} par "$\phi$ est monotone" et 
\eqref{eq160} par "$
\forall s\in [0,L],\quad 
 |\phi(s)|\in [0,|\Omega|]$". On cherche toujours la courbe qui minimise le maximum de la courbure géométrique.
\ifcase \cras
\end{remark}
\or
\end{remarque}
\fi

\ifcase \cras
\begin{remark}
\or
\begin{remarque}
\fi
\label{defiErem}
Notons que les résultats de Dubins sont valables pour tout couple de points $(A,B)$ et pour 
tout couple de vecteurs unitaires
$\left(\vec \alpha, \vec \beta\right)$. Ici, on impose les conditions supplémentaires \eqref{eq10tottot} et \eqref{eq20}.
On pourrait croire que le point $O$ peut être construit à partir des points distincts $A$ et $B$ et des vecteurs
$\vec \alpha$ et $\vec \beta$ 
deux vecteurs unitaires donnés vérifiant 
 \eqref{eq20} de la façon suivante :
\begin{equation}
\label{defiEremeq01}
\text{$O$ est l'unique intersection des droites passant respectivement par $A$ et $B$ et dirigées par 
$\vec \alpha$ et $\vec \beta$.}
\end{equation}
Si on le  définit ainsi, le point $O$ n'est pas nécessairement distinct de $A$ et de $B$ 
et \eqref{eq10tot} est alors remplacé \textit{a priori} par 
\begin{equation}
\label{eq10totbis}
\vec \alpha=\pm \frac{1}{OA}\overrightarrow{AO},\quad 
\vec \beta=\pm \frac{1}{OB}\overrightarrow{OB},\quad 
\end{equation}
Cette généralisation est inutile, comme le montre le lemme \ref{lemmeA}, qui servira à plusieurs reprises. 
\ifcase \cras
\end{remark}
\or
\end{remarque}
\fi
\ifcase \cras
\begin{lemma}
\or
\begin{lemme}
\fi
\label{lemmeA}
Soient deux vecteurs unitaires $\vec \alpha$ et $\vec \beta$, vérifiant \eqref{eq20} et  $O$ donné par 
\eqref{defiEremeq01}. 
S'il existe une courbe $X$ du plan 
vérifiant 
\eqref{eq10totnew},
\eqref{eq20},
\eqref{eq100tot},
\eqref{eq110},
\eqref{eq140} (ou \eqref{eq150}), 
alors $O$ est distinct de $A$ et de $B$ et si on considère les réels $u_0$ et $v_0$ tels que
$\overrightarrow{AO}=u_0 \vec \alpha$ et $\overrightarrow{OB}=v_0 \vec  \beta$, alors $u_0 $ et $v_0$  sont strictement positifs. 
\ifcase \cras
\end{lemma}
\or
\end{lemme}
\fi
\ifcase \cras
Voir la preuve du lemme \ref{lemmeA} en annexe \ref{preuve}, page \pageref{lemmeApreuve}.
\or
Voir la preuve du lemme \ref{lemmeA} en annexe \ref{preuve}.
\fi

\ifcase \faussepreuve

\ifcase \cras

\section{Existence d'une courbe minimisant  le maximum de la courbure}
\label{existencetheonewjuste}

On suppose désormais que sont fixés $O$ et $A$, $B$, deux à deux distincts  et $\vec \alpha$ et $\vec \beta$ vérifiant \eqref{eq10tottot} et \eqref{eq20}.

Donnons tout d'abord les deux estimations uniformes suivantes : 

\begin{proposition}
\label{existencetheoprop01new}
Il existe trois constantes $\mu$, $\nu$ et $\delta$ strictement positives 
 ne dépendant que des points $O$, $A$ et $B$, telles que, pour toute courbe $X$ de $\mathcal{E}$ : 
\begin{subequations}
\label{existencetheoprop01neweq}
\begin{align}
\label{lemmeBeq01}
&\mu \leq L(X) \leq \nu,\\
\label{existencetheolem02dfg}
&\vnorm[L^{\infty}(0,L)]{\vnorm{X''}}\geq \delta.
\end{align}
\end{subequations}
\end{proposition}

Voir preuve en annexe \ref{preuve}, page \pageref{preuveexistencetheoprop01new}.

On a alors le 
premier  résultat essentiel de cet article : 

\begin{theorem}
\label{existencetheoprop01juste}
L'ensemble $\mathcal{E}$ est non vide et il existe une courbe $X$ de $\mathcal{E}$ qui minimise le maximum de la courbure, c'est-à-dire vérifiant 
\eqref{eq200}.
\end{theorem}

Voir preuve en annexe \ref{preuve}, page \pageref{existencetheoprop01justepreuve}.

\fi

\or

\section{Existence d'une courbe minimisant  le maximum de la courbure}
\label{existencetheo}

\textbf{Attention,  ce résultat est faux !!!!! ; voir détails en annexe \ref{preuve}, page \pageref{refpreuvefausse}. Résultat finalement rétabli !!!!}

On suppose désormais que sont fixés $O$ et $A$, $B$, deux à deux distincts  et $\vec \alpha$ et $\vec \beta$ vérifiant \eqref{eq10tottot} et \eqref{eq20}.

On a alors le 
\ifcase \cras
 résultat essentiel : 
\or
premier  résultat essentiel de cette Note :
\fi

\ifcase \cras
\begin{theorem}
\or
\begin{theoreme}
\fi
\label{existencetheoprop01}
L'ensemble $\mathcal{E}$ est non vide et il existe une courbe $X$ de $\mathcal{E}$ qui minimise le maximum de la courbure, c'est-à-dire vérifiant 
\eqref{eq200}.
\ifcase \cras
\end{theorem}
\or
\end{theoreme}
\fi

Nous procédons en plusieurs étapes. 
Nous montrons tout d'abord que $\mathcal{E}$ est non vide. Ensuite, on montre qu'il existe 
$\mu$ et $\nu$, ne dépendant que des points $O$, $A$ et $B$, telles que 
\begin{equation}
\label{lemmeBeq01}
\forall X\in \mathcal{E},
\quad
\mu \leq L(X) \leq \nu.
\end{equation}
Puis, on montre qu'il  existe $\delta>0$ tel que,  pour toute courbe $X$ de
 $\mathcal{E}$,  
\begin{equation} 
\label{existencetheolem02dfg}
\vnorm[L^{\infty}(0,L)]{\vnorm{X''}}\geq \delta.
\end{equation}
On a donc
\begin{equation} 
\label{existencetheolem02dfgbis}
\delta_0=\inf_{X \in \mathcal{E}} \vnorm[L^{\infty}(0,L)]{\vnorm{X''}}>0.
\end{equation}
Enfin, on adapte la démonstration du résultat analogue de Dubins (\cite[Proposition 1]{MR0089457}) pour montrer que $\delta_0$ est atteint par une courbe $X$ de  $\mathcal{E}$
La différence essentielle est que la quantité à minimiser est la longueur de la courbe et non le maximum de la courbure et que les contraintes y sont moins restrictives.

Voir preuve en annexe \ref{preuve}.

\fi

\section{Construction de l'unique courbe de $\mathcal{E}$ formée d'un segment de droite et d'un arc de cercle.}
\label{uniquecercledroite}

\begin{definition}
\label{uniquecercledroitedef01}
On se donne  $O$ et $A$, $B$, deux à deux distincts  et $\vec \alpha$ et $\vec \beta$ vérifiant \eqref{eq10tottot} et \eqref{eq20}.
Nous dirons que nous sommes dans le cas symétrique si $OA =OB$ et dans le cas non symétrique si $OA\not =OB$.
\end{definition}

\ifcase \cras
\begin{lemma}
\or
\begin{lemme}
\fi
\label{uniquecercledroitelem01}
Il existe une unique courbe $X$ de $\mathcal{E}$
formée d'un arc de cercle de rayon $R_a>0$ et de longueur appartenant à  $]0,R_a\pi[$ 
dans le cas symétrique
et 
formée d'un arc de cercle de rayon $R_a>0$ et de longueur appartenant à  $]0,R_a\pi[$ 
et d'un segment de droite
de longueur non nulle dans le cas non symétrique. 
Le rayon $R_a$ du cercle est unique. Il ne dépend que de $O$, $A$ et $B$ et on a 
\begin{equation}
\label{uniquecercledroitelem01eq01}
\vnorm[L^{\infty}(0,L)]{\vnorm{X''}}=\frac{1}{R_a}.
\end{equation}
Pour toute la suite, cette courbe est notée sous la forme $X=\mathcal{J}(O,A,B)$ et le réel $R_a$ sous la forme $R_a(O,A,B)$.
\ifcase \cras
\end{lemma}
\or
\end{lemme}
\fi

La démonstration se fait de façon purement géométrique et est donnée 
\ifcase \cras
en annexe \ref{preuve}, page \pageref{uniquecercledroitelem01preuve}.
\or
en annexe \ref{preuve}.
\fi

\ifcase \cras
\begin{example}
\or
\begin{exemple}
\fi
\label{exemple02}
Traitons le cas particulier donné par 
$
A=(1/2,-1/2),\quad
O=(0,0),\quad
B=(0,-1/2),\quad
\linebreak[1]
\Omega=3\pi/4$.
\begin{figure}[h] 
\psfrag{O}{$O$}
\psfrag{A}{$A$}
\psfrag{B}{$B$}
\psfrag{C}{$C$}
\psfrag{D}{$D$}
\psfrag{a}{$\alpha$}
\psfrag{t}{$\theta$}
\begin{center} 
\epsfig{file=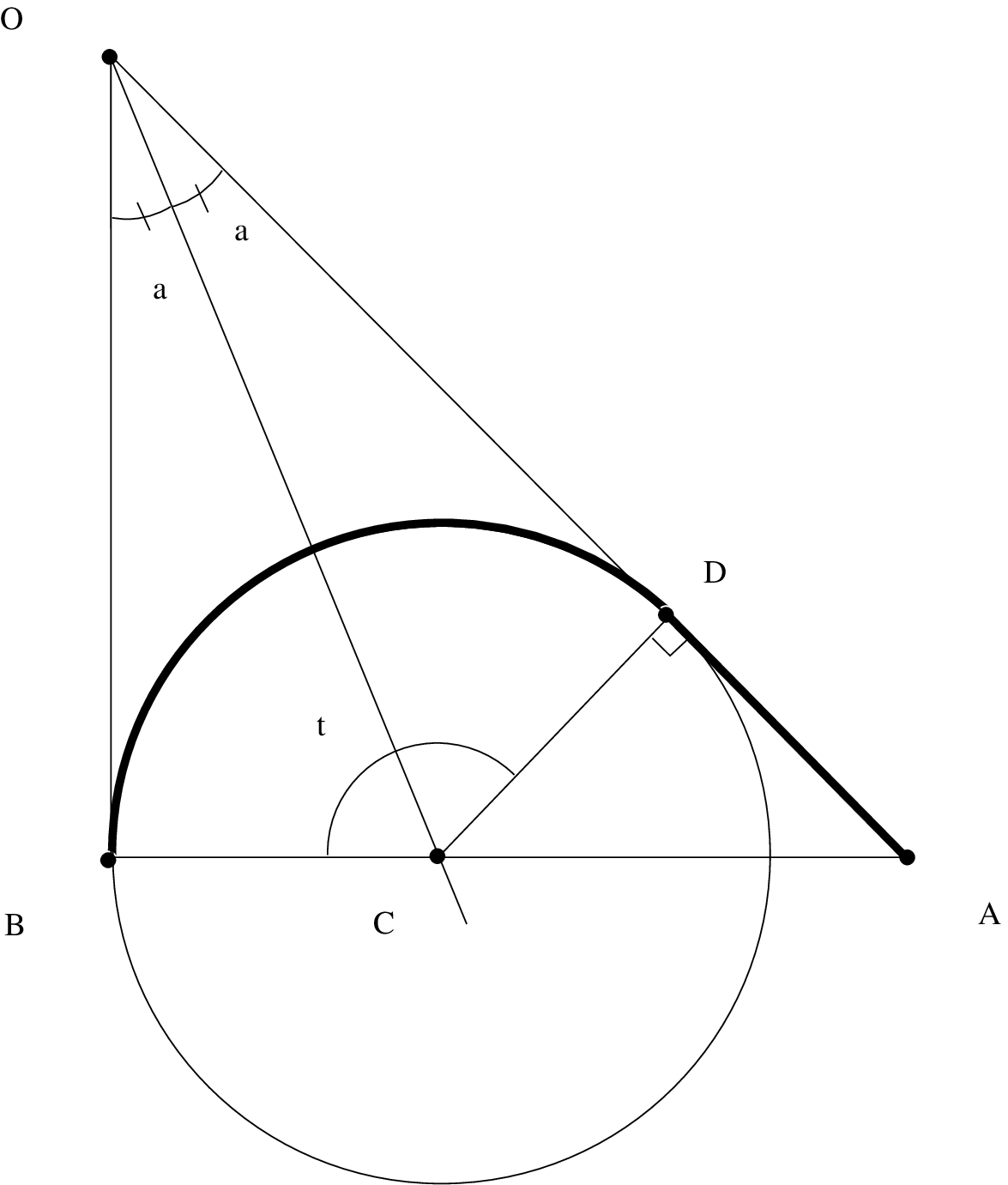, width=4  cm} 
\end{center} 
\caption{\label{forme5poptdet}La description géométrique de la courbe formée d'un arc de cercle et d'un segment de droite.\ifcase \cras \or /\textit{The geometric description of the curve composed of an arc of circle and a line segment.}\fi}
\end{figure}
La construction reprend la méthode donnée dans la preuve du lemme \ref{uniquecercledroitelem01}. 

Comme indiqué sur la figure \ref{forme5poptdet}, la courbe constituée par  la réunion  d'un arc de cercle et d'un segment de droite est définie de la façon suivante (le triangle $OAB$ étant  isocèle rectangle en $B$ avec $OB=BA=1/2$) :  
$(OC)$ est la bissectrice de l'angle $\widehat O$ avec $\alpha=\pi/8$ ;
$\theta =3\pi/4$ ; 
l'arc de cercle a pour centre  $C$ et pour rayon $R$ donné par 
$
R=\frac{\sqrt{2}-1}{2}\approx  0.20710678,
$
 et est limité par les points $B$  et $D$ ;
le segment de droite est le segment $[DA]$ avec $DA=(\sqrt{2}-1)/2$.

La parabole définie dans le cadre du brevet \cite{brevetJB} a été remplacée par cette courbe et cela nous a permis  de faire passer le minimum du rayon de courbure
de $1/25\,\sqrt {5} \approx 0.089$ à $(\sqrt{2}-1)/2\approx 0.207$.
\ifcase \cras
\end{example}
\or
\end{exemple}
\fi

En reprenant la notation \eqref{rapeldubinseq01}, construisons maintenant autrement la courbe de $\mathcal{E}$, définie dans le lemme \ref{uniquecercledroitelem01}, en utilisant
les courbes de Dubins.

\ifcase \cras
\begin{lemma}
\or
\begin{lemme}
\fi
\label{optdubinslem01}
Considérons 
$
\mathcal{F}=
\left\{R\in \Erps, \quad   \mathcal{G}(R)\in     \mathcal{E}\right\}
$.
Si $R_a(O,A,B)$ est le nombre défini dans le lemme \ref{uniquecercledroitelem01}, alors 
$
\mathcal{F}=]0,R_a(O,A,B)]
$
et la courbe de Dubins   
$\mathcal{G}(R_a(O,A,B))$  
est l'unique courbe définie dans le lemme \ref{uniquecercledroitelem01}.
\ifcase \cras
\end{lemma}
\or
\end{lemme}
\fi

Autrement dit, $\mathcal{G}(R_a(O,A,B))$ est optimale : elle correspond à la plus grande valeur possible de $R$, pour laquelle 
la courbe de Dubins
$\mathcal{G}(R)$ est dans $\mathcal{E}$.
Voir la figure \ref{coducl1} qui correspond au cas optimal.
La démonstration se fait de façon purement géométrique et est donnée en 
\ifcase \cras
annexe \ref{preuve}, page \pageref{optdubinslem01preuve}.
\or
annexe \ref{preuve}.
\fi

\ifcase \cras
\or
\input{vrai_probleme_pointB}
\fi

\ifcase \cras

\section{Présentation du problème discret}
\label{discret}

Les techniques utilisées par Dubins dans \cite{MR0089457} ne peuvent être utilisées pour déterminer de façon explicite 
 l'une des courbes optimales du théorème 
 ; 
en effet, on ne minimise pas la longueur de la courbe et surtout, le rayon de courbure maximal, 
connu chez Dubins, nous est inconnu. L'idée de cette section est de remplacer le problème initial par un problème discret, où il est possible de déterminer explicitement 
le rayon de courbure maximal et de déterminer les paramètres optimaux pour maximiser  ce  rayon de courbure minimal.

On considère le repère orthonormé $\left(A, \vec {\alpha}, \vec {k}\right)$, et en notant $X(s)=(  x(s),  y(s))$, 
les coordonnées des points de  la courbe $X$ dans ce repère, 
on a les relations habituelles  
\begin{subequations}
\label{eqderphinbbbb}
\begin{align}
\label{eqderphinab}
&\frac{d  x}{ds}=\cos   \phi,\\
\label{eqderphinbb}
&\frac{d y}{ds}=\sin   \phi.
\end{align}
\end{subequations}

Remarquons tout d'abord que, définir $X\in \mathcal{E}$
revient donc à définir successivement 
\begin{itemize}
\item[$\bullet$]
$L>0$ ; 
\item[$\bullet$]
la courbure $c\in L^\infty(0,L)$, presque partout positive ; 
\item[$\bullet$]
l'angle $\phi\in W^{1,\infty}(0,L)$ grâce à 
\begin{equation}
\label{passagecontinuediscreteq01}
\forall s\in [0,L],\quad
\phi(s)=\int_0^s c(u) du,
\end{equation}
qui doit vérifier 
\begin{equation}
\label{passagecontinuediscreteq05}
\int_0^L c(u) du=\Omega,
\end{equation}
\item[$\bullet$]
les fonctions $x$ et $y$ de  $W^{2,\infty}(0,L)$ 
définies respectivement par
\begin{equation}
\label{passagecontinuediscreteq10}
\forall s\in [0,L],\quad
x(s)=\int_0^s \cos\phi(u) du,\quad
y(s)=\int_0^s \sin\phi(u) du,
\end{equation}
qui doivent vérifier
\begin{equation}
\label{passagecontinuediscreteq20}
\int_0^L \cos\phi(s) du=x_B,\quad 
\int_0^L \sin\phi(s) du=y_B.
\end{equation}
\end{itemize}

Une discrétisation naturelle consiste à définir alors $c$ constante par morceaux, sur un nombre fini d'intervalles $[L_i,L_{i+1}]$.
Supposons que $p$ désigne le nombre d'intervalles et que, pour tout $i\in \{0,...,p-1\}$, $c$ est constant et égal à $1/R_i$ avec $R_i>0$. 
Les inconnues sont alors $\textit{a priori}$, $p\in \En^*$, $L>0$, $R=(R_0,...,R_{p-1})\in {\left(\Erps\cup \{+\infty\}\right)}^p$  et 
$L=(L_0,...,L_{p-1},L_p)\in {\Erps}^{p+1}$,  avec $L_0=0$, qui définissent des fonctions qui doivent vérifier 
\eqref{passagecontinuediscreteq05} et \eqref{passagecontinuediscreteq20}.
La courbe $X$ est donc constituée d'un nombre fini d'arcs de cercles ou de segments de droite.

La fonction $\phi$, définie par \eqref{passagecontinuediscreteq01}, est linéaire par morceaux, continue sur $[0,L]$. Si on pose 
\begin{equation}
\label{passagecontinuediscreteq30}
\forall i\in \{0,...,p\},\quad 
\phi_i=\phi(L_i),
\end{equation}
on a, d'après \eqref{passagecontinuediscreteq01}, 
\begin{equation}
\label{passagecontinuediscreteq40}
\forall i\in \{0,...,p-1\},\quad 
\phi_{i+1}-\phi_i=\frac{L_{i+1}-L_i}{R_i}.
\end{equation}
Ce problème discret doit
vérifier les deux contraintes \eqref{passagecontinuediscreteq05} et \eqref{passagecontinuediscreteq20}
ce qui donne lieu à un problème \textit{a priori} non linéaire. 
Cette
discrétisation est décrite en section \ref{discretgen}.
Ensuite, 
pour simplifier, 
nous considérons  qu'un intervalle où $c$ est nul (correspondant à 
un segment de droite constituant une partie de $X$) peut être remplacé par un intervalle où $c$ est "petit" et strictement positive, 
le segment étant approché par un cercle de "grand" rayon. 
Enfin, nous supposons $\phi_{i+1}-\phi_i=\xi$ constant. On a donc $\phi_{p}=\Omega =p\xi$ et donc $\xi=\Omega/p$. 
Enfin, \eqref{passagecontinuediscreteq40} donne 
$L_{i+1}-L_i=R_i \xi$. 
Il est ensuite plus simple de considérer l'affixe de $b$ de  $B$ et de constater que la contrainte 
\eqref{passagecontinuediscreteq20}
 fournit un problème 
 linéaire. 
Cette
discrétisation est décrite en section \ref{discretcaspthencon}.

\subsection{Cas général}
\label{discretgen}

On considère un entier $p\in \En^*$
et on construit  $p+1$ points du plan, notés ${(A_k)}_{0\leq k\leq p}$,
reliés par $p$ courbes ${\mathcal{D}}_k$ qui seront soit des arc de cercle, soit des segments de droites.
Notons, pour $0\leq k \leq p$, $a_k$ l'affixe complexe du point $A_k$,  et $L_k>0$, la longueur de la courbe ${\mathcal{D}}_k$.

\begin{proposition}
\label{discretprop01}
Notons  $b$ l'affixe respective de $B$  et  $\beta$ l'affixes du vecteur  $\vec \beta$. Soit $p\in \En^*$. Pour tout ${(\theta_k)}_{0\leq k\leq p-1}\in {\left(\Erp\right)}^{p+1}$, on pose 
\begin{equation*}
\forall k\in \{0,...,p\},\quad
\phi_k=\sum_{j=0} ^{k-1}
\theta_j,
\end{equation*}
avec pour convention $\phi_0=0$ et 
\begin{equation*}
\forall k\in \{0,...,p\},\quad
\gamma_k=
\frac{\sin(\theta_k/2)}{\theta_k/2}e^{i\theta_k/2},
\end{equation*}
prolongé par continuité par convention par $1$ si $\theta_k=0$. 
La courbe formée par la réunion des 
$p$ courbes ${\mathcal{D}}_k$ 
définies plus haut,
par les $2p$ réels 
${(\theta_k)}_{0\leq k\leq p-1}\in {(\Erp)}^{p}$ et ${(L_k)}_{0\leq k\leq p-1}\in {(\Erps)}^{p}$
appartient à $\mathcal{E}$ ssi 
\begin{subequations}
\label{discretprop01eq20tot}
\begin{align}
\label{discretprop01eq20a}
&\phi_p=\Omega,\\
\label{discretprop01eq20b}
&\forall k\in \{0,...,p\},\quad 
u_k=e^{i\phi_k},\\
\label{discretprop01eq20bc}
&b=\sum_{k=0}^{p-1} \gamma_ku_kL_k.
\end{align}
\end{subequations}
\end{proposition}

Il suffit de remarquer que
la courbe formée par la réunion des 
$p$ courbes ${\mathcal{D}}_k$  est dans $\mathcal{E}$ ssi
$u_0=1$ et  $u_{p}=\beta$, $A_0=A$, $A_p=B$.
Notons que la courbure de cette courbe est définie partout (sauf aux points $A_k$) et est constante par morceaux, valant $1/R_k$. Ainsi, 
\eqref{eq140} est vérifiée.

\subsection{Cas $p$ quelconque et $\theta_k$ constant}
\label{discretcaspthencon}

\begin{proposition}
\label{discretcaspthenconprop01}
Soit $p\in \En^*$. On définit ${(\theta_k)}_{0\leq k \leq p-1}$ par 
\begin{equation}
\label{discretcaspthenconprop01eq01}
\forall k\in \{0,...,p-1\},\quad
\theta_k=\frac{\Omega}{p}\in \Erps.
\end{equation}
Dans ce cas, la courbe $X$ définie dans 
la proposition \ref{discretprop01}
est formée par $p$ arcs de cercles de rayons $R_k>0$, pour $0\leq k\leq p-1$. 
Elle appartient à $\mathcal{E}$ ssi 
\begin{equation*}
b=\sum_{k=0}^{p-1} 
i\left(1-e^{i\theta_0}\right)e^{ik\theta_0}  R_k,
\end{equation*}
ce qui est équivalent au système linéaire
\begin{equation}
\label{discretcaspthenconprop01eq20}
\mathcal{A}R=\mathcal{B},
\end{equation}
où
\begin{subequations}
\label{discretcaspthenconprop01eq30tot}
\begin{align}
\label{discretcaspthenconprop01eq30a}
&\mathcal{A}\in {\mathcal{M}}_{2,p}\left(\Er\right),\\
\label{discretcaspthenconprop01eq30b}
\text{avec},  \forall k\in \{1,...,p\}, \quad & \mathcal{A}_{1,k}=\text{Re}\left(i\left(1-e^{i\theta_0}\right)e^{i(k-1)\theta_0} \right),\quad 
\mathcal{A}_{2,k}=\text{Im}\left(i\left(1-e^{i\theta_0}\right)e^{i(k-1)\theta_0} \right),\\
\label{discretcaspthenconprop01eq30d}
&\mathcal{B}\in {\mathcal{M}}_{2,1}\left(\Er\right),\\
\label{discretcaspthenconprop01eq30e}
\text{avec}\quad & \mathcal{B}_{1}=\text{Re}(b),\quad 
\mathcal{B}_{2}=\text{Im}(b),\\
\label{discretcaspthenconprop01eq30g}
&R=(R_0,...,R_{p-1})\in {\left(\Erps\right)}^{p}.
\end{align}
\end{subequations}
Enfin, on a 
\begin{equation}
\label{discretcaspthenconprop01eq40}
\vnorm[L^{\infty}(0,L)]{\vnorm{X''}}=
\frac{1}{\displaystyle{\min_{0\leq k\leq p-1} R_k}}.
\end{equation}
\end{proposition}

Avec l'identification de ${\mathcal{E}}_p$  à une partie de ${\left(\Erps\right)}^{p}$, on a 
\begin{equation*}
\vnorm[L^{\infty}(0,L)]{\vnorm{X''}}=
\max_{0\leq k\leq p-1} \frac{1}{|R_k|}=\max_{0\leq k\leq p-1} \frac{1}{R_k},
\end{equation*}
et dont découle 
\eqref{discretcaspthenconprop01eq40}.
Ainsi minimiser $\vnorm[L^{\infty}(0,L)]{\vnorm{X''}}$ revient à maximiser $\displaystyle{\min_{0\leq k\leq p-1} R_k}$.

Il suffit ensuite d'appliquer la proposition \ref{discretprop01} avec 
\eqref{discretcaspthenconprop01eq01}, dans laquelle on  a, pour tout $k$, $\theta_k=\theta_0$ et $\phi_k=k\theta_0$

On peut alors, à $p$ fixé, adopter la définition suivante : 

\begin{definition}
\label{discretcaspthencondefi01}
On appelle  ${\mathcal{E}}_p$ l'ensemble des courbes du plan formées de 
 $p$ arcs de cercles de rayons $R_k>0$, pour $0\leq k\leq p-1$, comme défini dans 
les propositions 
\ref{discretprop01}
et 
\ref{discretcaspthenconprop01}. Chaque courbe $X$ de cet ensemble pourra être identifiée
au vecteur $R$ de ${\left(\Erps\right)}^{p}$ défini par \eqref{discretcaspthenconprop01eq30g} et l'ensemble 
${\mathcal{E}}_p$ pourra donc être identifié à une partie de ${\left(\Erps\right)}^{p}$.
\end{definition}

On a alors le 
second   résultat essentiel de cet article :

\begin{theorem}
\label{existencetheoprop01discprop01}
Soit $p\geq 2$.
Si ${\mathcal{E}}_p$ est non vide,
il existe une courbe $X$ de ${\mathcal{E}}_p$ qui minimise le maximum de la courbure, c'est-à-dire vérifiant 
\begin{equation}
\label{existencetheoprop01eqdisceq01}
\vnorm[L^{\infty}(0,L)]{\vnorm{X''}}=
\min_{Z\in {\mathcal{E}}_p} \vnorm[L^{\infty}(0,L(Z))]{\vnorm{Z''}}.
\end{equation}
Pour $X=(X_1,...,X_p) \in\Erp^p$, on note 
\begin{equation}
\label{existencetheoprop01eqdisceq10}
F(X)=\min_{1\leq k\leq p} X_k. 
\end{equation}
En identifiant ${\mathcal{E}}_p$  à une partie de ${\left(\Erps\right)}^{p}$,
le problème \eqref{existencetheoprop01eqdisceq01} est aussi équivalent à 
\begin{equation}
\label{existencetheoprop01eqdisc02}
F(X)=
\max_{\substack{Z\in {\left(\Erps\right)}^{p},\\ \mathcal{A}Z=\mathcal{B}}} F(Z),
\end{equation}
qui admet une solution, 
s'il existe $X^0\in {\left(\Erps\right)}^{p}$ tel que $\mathcal{A}X^0=\mathcal{B}$.
\end{theorem}

Voir preuve en annexe \ref{preuve}, page 
\pageref{existencetheoprop01discprop01preuve}.
\ifcase \faussepreuve
Il s'agit simplement d'utiliser la proposition \ref{existencetheoprop01new}. 
\fi

Il faut supposer ${\mathcal{E}}_p$ non vide,
ce que l'on ne peut pas montrer en théorie, mais que l'on pourra vérifier numériquement.
Pour résoudre \eqref{existencetheoprop01eqdisc02},
nous utiliserons par exemple la fonction \path|fminimax| de 
 Matlab avec l'option 'Medium-Scale'.
Cet algorithme utilise la  méthode 'sequential quadratic programming' (SQP) présentée dans \cite{Brayton1979}.
La solution numérique déterminée dans cette proposition sera notée $X^d$.

On peut aussi déterminer la courbe de Dubins discrète obtenue en minimisant cette fois-ci la longueur
discrète et en imposant que $\min(X)\geq \min(X^d)$ cela grâce à la fonction \path|linprog| de matlab, 
fondée sur la méthode 'Linear Interior Point Solver' (LIPSOL), présentée dans \cite{MR1671584}.
La solution obtenue est notée $X^c$.

\subsubsection{Simulations numériques}\
\label{discretsimul}

\newcommand{\taillefig}{7}

Reprenons  l'exemple \ref{exemple02}  pour 
$p=300$.

\begin{figure}
\psfrag{Courbes}{}
\psfrag{Rayons de courbure}{}
\centering
\subfigure[\label{exemple12021}Courbes]%
{\epsfig{file=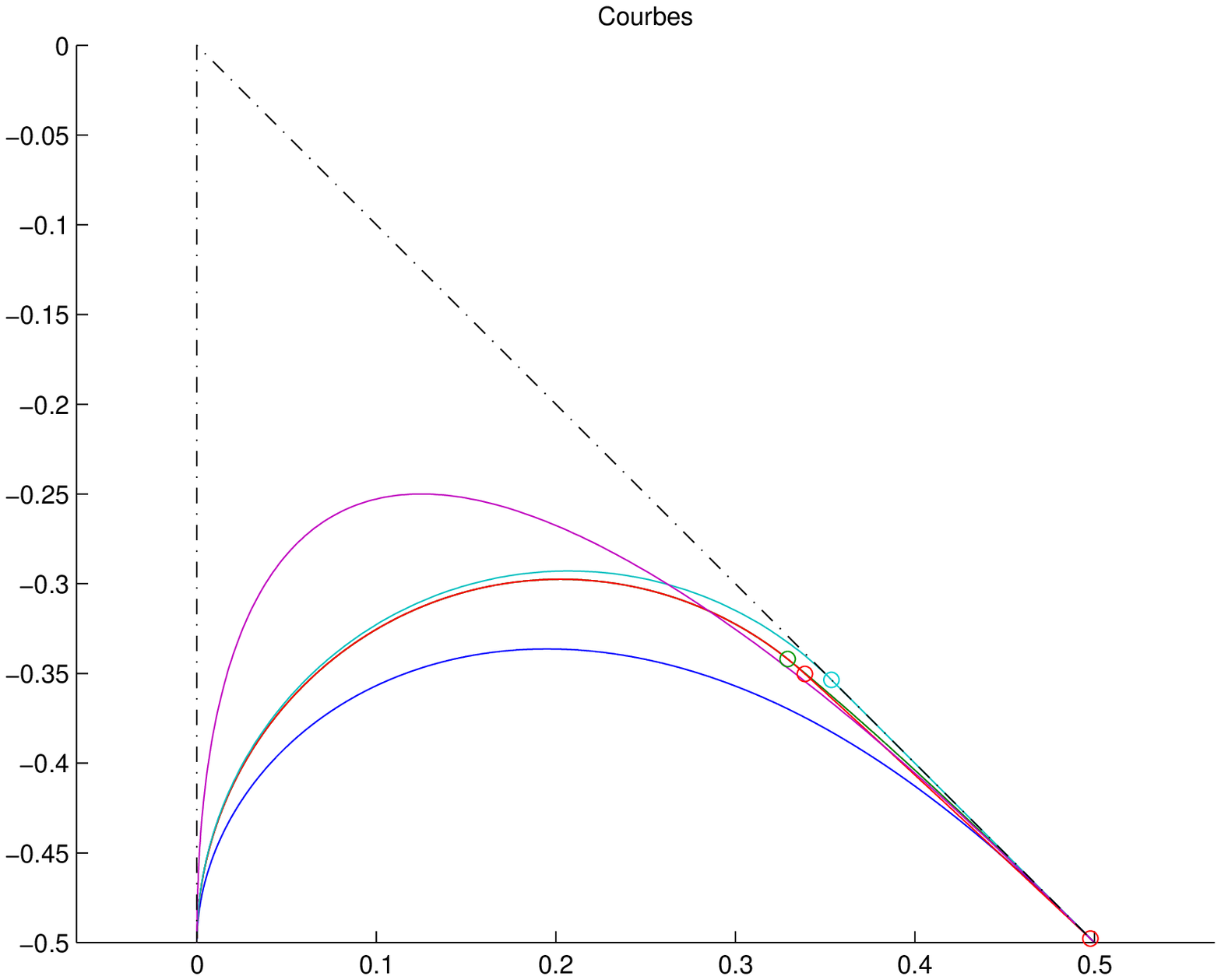, width=\taillefig cm}}
\subfigure[\label{exemple12020}Rayons de courbure]%
{\epsfig{file=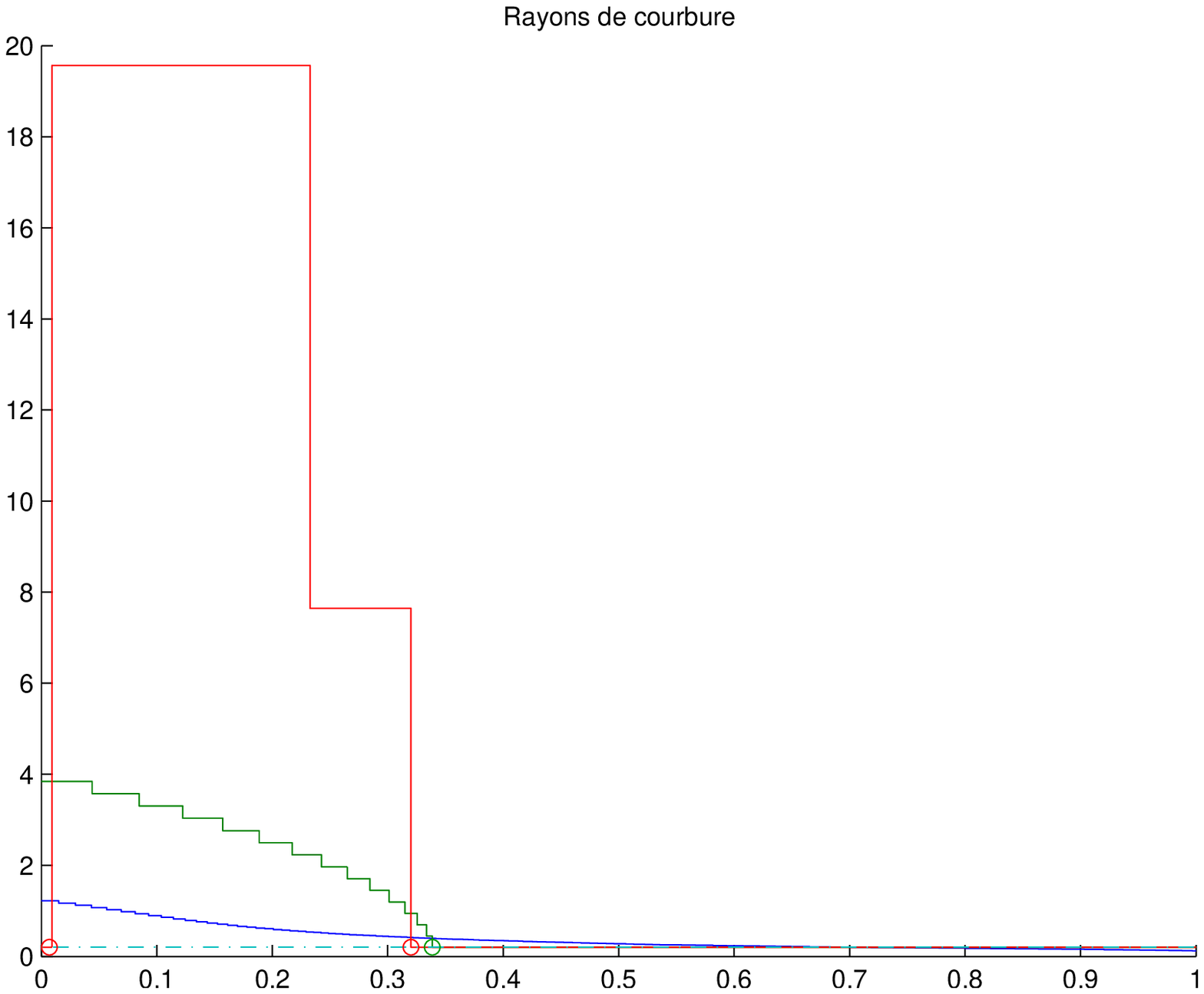, width=\taillefig cm}}
\caption{\label{exemple120tot}Les courbes obtenues.
Sur  figure \ref{exemple12021}, 
on a représenté $X^0$ en bleu, $X^d$  en vert, $X^c$ en rouge, la courbe théorique du lemme \ref{uniquecercledroitelem01} en bleu clair 
et la parabole  en violet.
Sur la figure \ref{exemple12020} sont représentés les rayons de courbures  (constants par morceaux), avec les mêmes conventions de couleurs,
hormis la parabole. 
De petits arcs de cercles délimitent la partie circulaire et la partie quasirectiligne.%
\ifcase \cras \or //\textit{%
The  obtained curves.
The curves obtained for different values of $ p $. In the figure \ref{exemple12021}, the points $ A_k $ have been represented and the curves $ X^d $ and $ X^c $ are mingled.
We have plotted $ X^0 $ in blue, $ X^d $ in green, $ X^c $ in red, the theoretical curve of the lemma \ref{uniquecercledroitelem01} in light blue and the parabola in purple.
In the figure \ref{exemple12020}, the radius of curvature (constant piecewise) obtained are drawn with the same color convention, except the parabola.
Small arcs of circles delimit the circular part and the quasirectiligne part.%
}\fi}
\end{figure}

Voir 
la figure \ref{exemple120tot}, où ont été choisis les paramètres de l'exemple \ref{exemple02}. 
Cette figure met en évidence le fait que les courbes obtenues sont constituées de cercles de grands rayons de courbures, 
correspondant à la partie rectiligne de la courbe théorique, puis de cercles de rayons très proches du rayon de la partie circulaire.
La courbe présente semble proche donc de la courbe décrite en théorie.
en section \ref{uniquecercledroite}.

\fi

\ifcase \cras

\section{Construction effective de la pièce du circuit et exemple d'un circuit}
\label{construceffct}

Si on choisit les dimensions de la section standard Brio, choisis pour les rails Easyloop, on obtient donc finalement la 
pièce 6 représentée sur la figure \ref{fin}.

\begin{figure}
\centering
\subfigure[\label{fin} La forme optimale./\textit{The optimal shape.}] 
{\epsfig{file=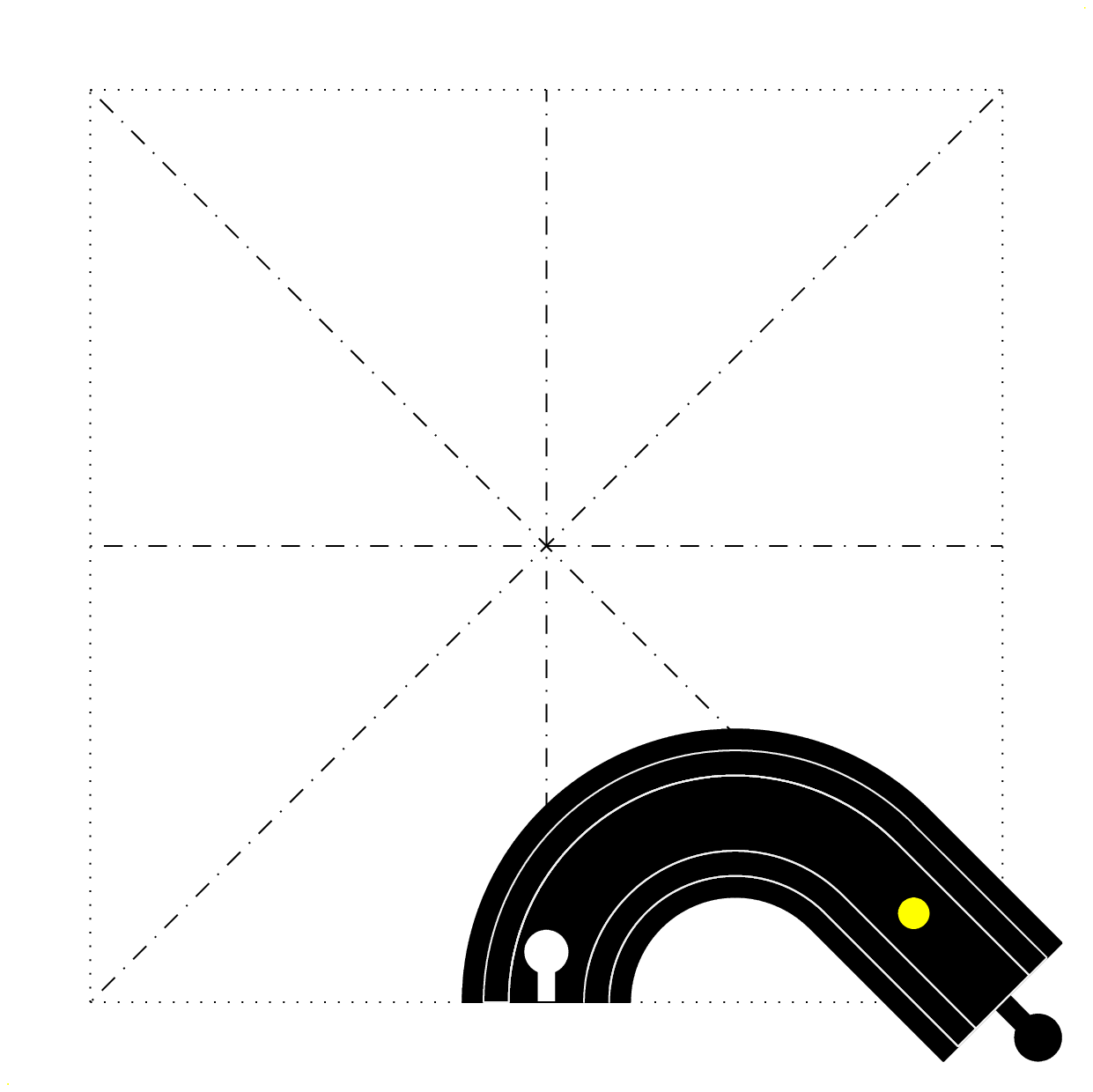, width=6 cm}}
\qquad
\subfigure[\label{exemplemcarder01b}Un exemple de circuit avec la pièce 6 optimale./\textit{An exemple of track with optimal shape 6.}]
{\epsfig{file=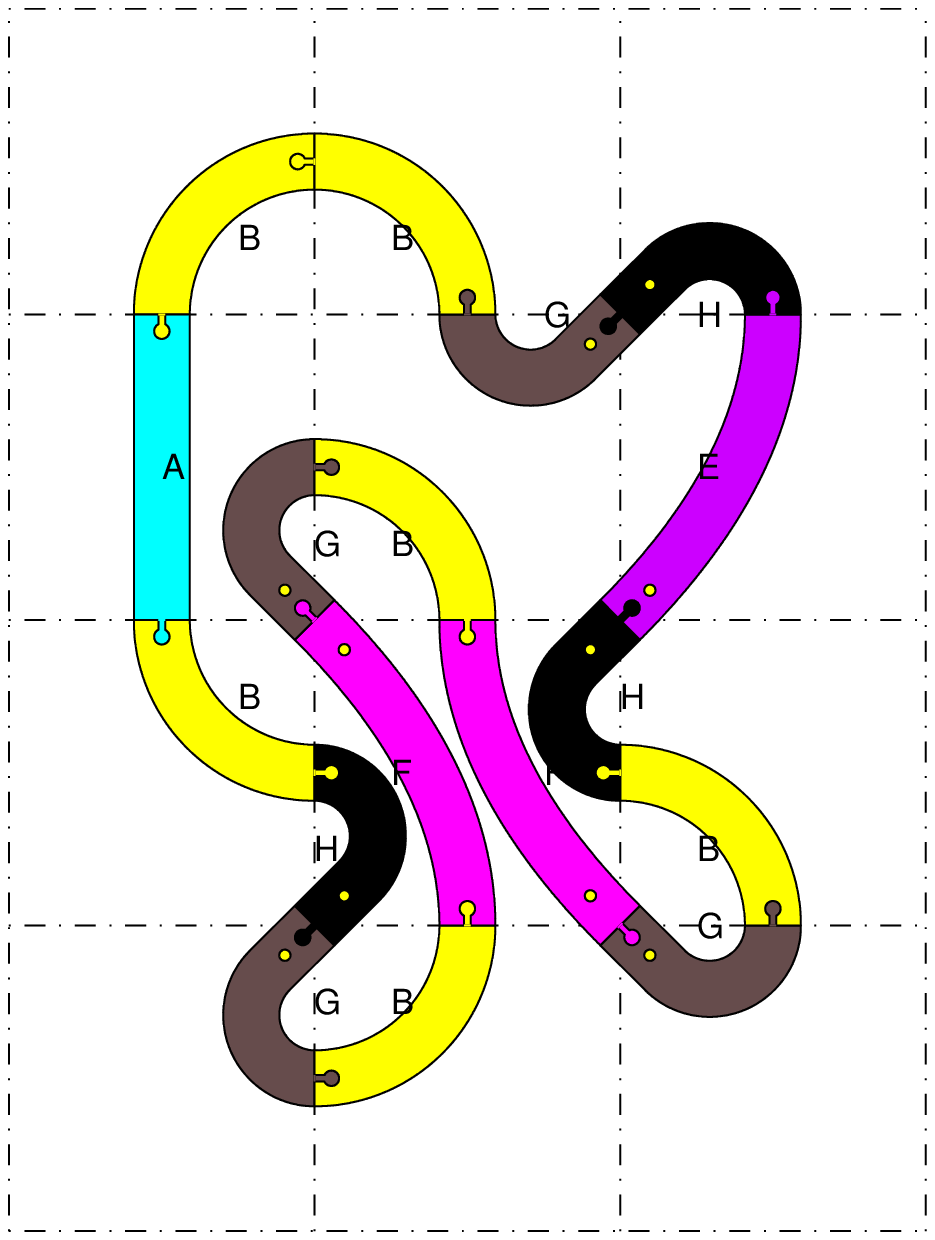, width=5 cm} }
\caption{\label{finexemplemcarder01b} Utilisations de la courbe optimale/\textit{Usings the optimal curve}}
\end{figure}

On pourra aussi consulter la figure \ref{exemplemcarder01b} qui montre un exemple d'un ciruit contenant la pièce optimale définie plus haut.

\fi

\section{Conclusion}
\label{conclusion}

\ifcase \cras
On a montré 
l'existence d'une courbe 
de 
\definitE\
qui maximise le minimum du rayon de courbure.
De plus, il 
existe une unique courbe ${\mathcal{J}}$ de $\mathcal{E}$,
 formée d'un segment de droite et d'un arc de cercle.
En étudiant les courbes constituées d'un nombre finis d'arcs de cercles, on a constaté numériquement que cette dernière courbe semble être ${\mathcal{J}}$.
\fi

\ifcase \cras
\or
\input{vrai_probleme_pointC}
\fi

\ifcase \cras
L'importante lacune de  
\ifcase \cras
cet article,
\or
cette Note, 
\fi
qu'il conviendrait de combler  par la suite, est de montrer 
\ifcase \cras
ce que l'on a observé numériquement, c'est-à-dire que 
la courbe
mise en évidence dans le théorème
\ifcase \faussepreuve
\ref{existencetheoprop01juste}
\or
\ref{existencetheoprop01}
\fi
est égale à celle construite dans le lemme \ref{uniquecercledroitelem01}.
Par la suite, il serait intéressant de montrer les éventuelles unicités des courbes, dans le cas continu comme dans les différents cas discrets évoqués.
\fi
\fi

\ifcase \cras
\or
\input{vrai_probleme_pointG}
\fi

\ifcase \cras
\or

\section*{Remerciements}

Je remercie vivement les rapporteurs anonymes pour la qualité et la précision de leurs remarques qui m'ont aidé à améliorer cette Note et les travaux préliminaires \cite{piece6_optimale_JB_2019_X4}.

\fi

\ifcase \cras

\printbibliography

\or

\bibliographystyle{unsrt}
\bibliography{piece6_optimale_JB_CRAS_2019_bis}

\fi

\appendix

\clearpage

\ifcase \cras
\markboth{J\'ER\^OME BASTIEN}%
{\ftitrefacul}
\fi

\ifcase \cras
\or
\textit{Les annexes \ref{preuve} 
 et \ref{preuvecomplete} sont destinées aux rapporteurs et ne font pas partie du projet de Notes.
Il est en effet  écrit sur le site des CRAS : 
"Il n'est parfois pas possible, en raison de la concision exigée, de donner les démonstrations ou les preuves complètes du résultat énoncé, spécialement dans la série I ; il est alors recommandé de joindre à l'appui de la Note un texte, si possible dactylographié, explicitant les compléments nécessaires à une bonne compréhension qui facilite l'examen de la Note par le présentateur et qui sera conservé cinq ans dans les Archives de l'Académie pour pouvoir être communiqué à tout lecteur des Comptes rendus qui en fait la demande. Il est alors fait mention au bas de la Note de l'existence de ce document par la formule suivante ; "Résumé d'un texte qui sera conservé cinq ans dans les Archives de l'Académie et dont copie peut être obtenue."  "
Les preuves de l'annexe \ref{preuve} viennent d'être déposées\ dans \url{https://arxiv.org}. Voir \cite{piece6_optimale_JB_2019_X4}.
L'annexe \ref{preuvecomplete} contient toutes les preuves des résultats essentiels.
Le cas échéant, ces preuves feront l'objet d'une publication ultérieure.}
\fi

\section{Preuves annexes}
\label{preuve}


\ifcase \cras
\or
La totalité des preuves de cette annexe est publiée dans \cite{piece6_optimale_JB_2019_X4}.
\fi

Rappelons tout d'abord que, comme les courbes représentatives  d'applications convexes, 
on a le lemme suivant :

\ifcase \cras
\begin{lemma}
\or
\begin{lemme}
\fi
\label{lemmeconvexite}
Soit une courbe $X$ vérifiant 
\eqref{eq10totnew},
\eqref{eq20},
\eqref{eq100tot} ,
\eqref{eq110} et  
\eqref{eq150}
alors, 
pour tout $s\in [0,L]$, la courbe est incluse dans le demi-plan délimitée par la droite tangente à la courbe au point $X(s)$, 
du côté de $\vec N(s)$,  la normale extérieure à la courbe en $X(s)$. 
\ifcase \cras
\end{lemma}
\or
\end{lemme}
\fi

\ifcase \cras

\ifcase \cras
\begin{proof}\
\or
\begin{preuve}\
\fi

\begin{figure}[h]    
\psfrag{A}{$A$}
\psfrag{B}{$B$}
\psfrag{a}{$\vec \alpha$}
\psfrag{b}{$\vec \beta$}
\begin{center} 
\epsfig{file=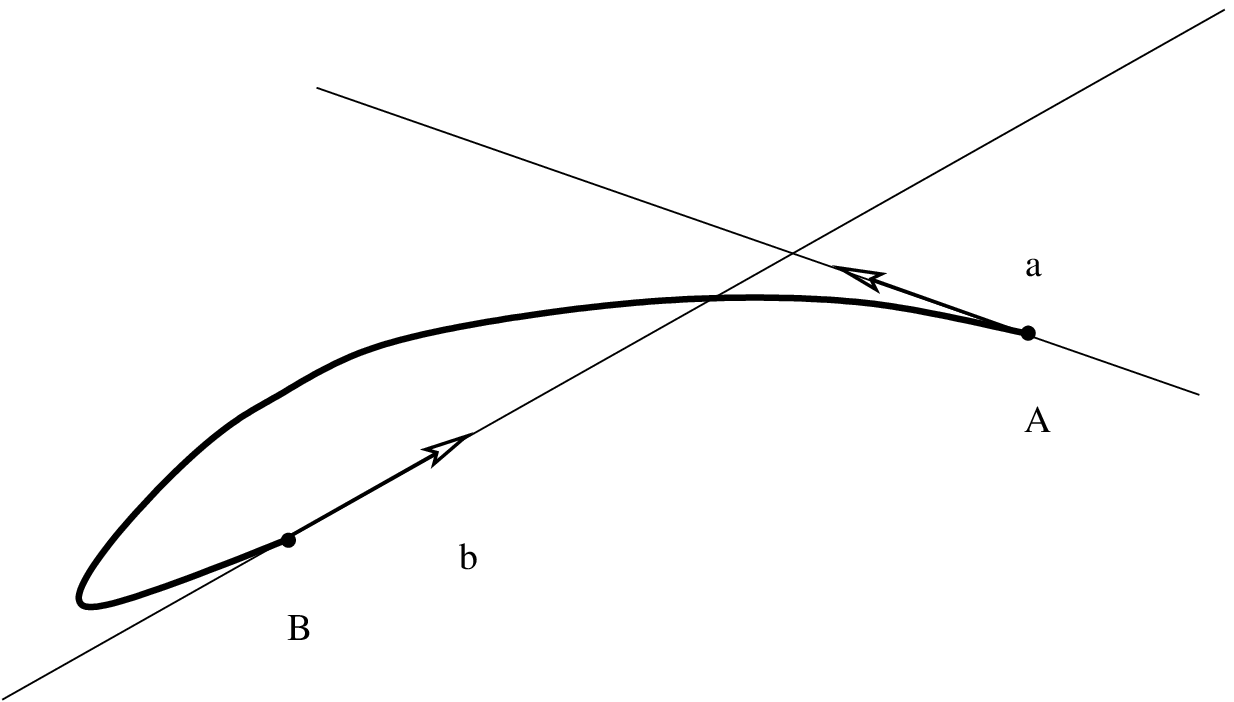, width=6 cm}  
\end{center} 
\caption{\label{contre_exemple_convexite}Contre-exemple au lemme  \ref{lemmeconvexite} sans l'hypothèse \eqref{eq20}.\ifcase \cras \or /\textit{Counter-example to Lemma  \ref{lemmeconvexite} without assumption 
\eqref{eq20}.}\fi}
\end{figure}

Notons que, sans l'hypothèse \eqref{eq20}, cette propriété devient fausse comme le montre le contre-exemple de la figure 
\ref{contre_exemple_convexite}.

\begin{figure}[h]    
\psfrag{a}{$X'(s)$}
\psfrag{b}[][l]{$\vec N(s)$} 
\psfrag{P}{$X(s)$}
\psfrag{Q}{$X(t)$}
\begin{center} 
\epsfig{file=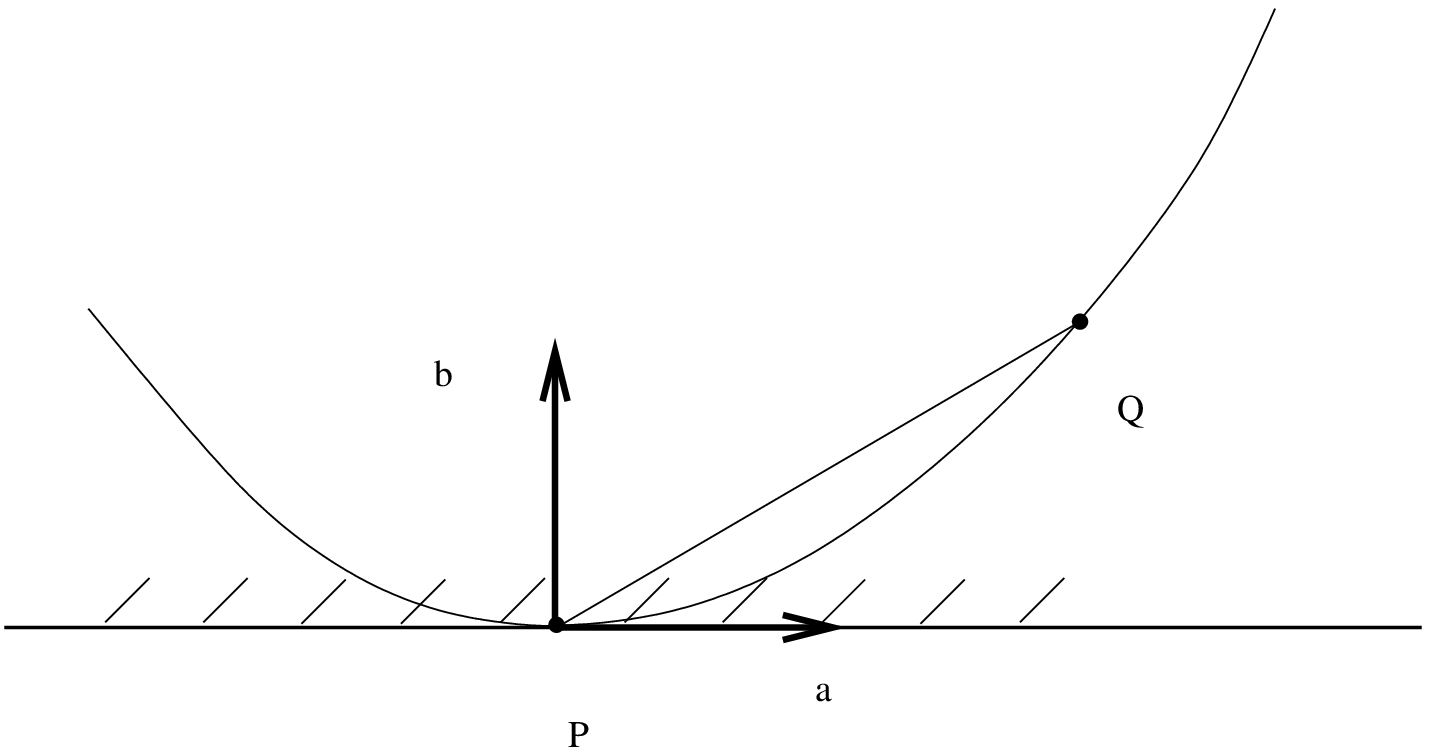, width=9 cm}  
\end{center} 
\caption{\label{courbe_meme_cote}La courbe est toujours du côté de la normale extérieure.\ifcase \cras \or /\textit{The curve is ever in the side of the outer-pointing normal.}\fi}
\end{figure}

Notons que 
\begin{equation*}
N(s)=\sigma(X'(s)),
\end{equation*}
où 
\ifcase \cras
$\sigma$    est la rotation vectorielle d'angle $\pi/2$.
\or
$\sigma$    est la rotation vectorielle d'angle $\pi/2$, déjà utilisée dans la preuve de la proposition \ref{lemme100}.
\fi
Voir sur la figure \ref{courbe_meme_cote}, la situation représentée. 
\ifcase \cras
On note $\prodsca{.}{.}$, le produit scalaire Euclidien de $\Er^2$ (qui induit la norme Euclidienne  $\vnorm{.}$  de $\Er^2$). 
\or
$\prodsca{.}{.}$ est le produit scalaire Euclidien de $\Er^2$.
\fi
Pour $s\in [0,L]$ fixé et pour tout $t\in [0,L]$, on pose 
\begin{equation*}
\gamma(t)=\prodsca{X(t)-X(s)}{\sigma(X'(s))},
\end{equation*}
dont la dérivée vaut d'après \eqref{eqderphinbbbb}
\begin{align*}
\gamma(t)
&=\prodsca{X'(t)}{\sigma(X'(s))},\\
&=
\begin{pmatrix}
\cos   \phi(t) \\ \sin   \phi(t)
\end{pmatrix}
.
\begin{pmatrix}
-\sin\phi(s)\\ \cos\phi(s)
\end{pmatrix},\\
&=
\sin   \phi(t) \cos\phi(s)-\cos   \phi(t)\sin\phi(s),
\end{align*}
et donc 
\begin{equation}
\label{ffvsjhgfss}
\forall t\in [0,L],\quad
\gamma'(t)=
\sin\left(\phi(t)-\phi(s)\right).
\end{equation}
Or,  d'après 
\eqref{eq20},
 et  
\eqref{eq150},
on a, pour tout $t$, 
\begin{equation*}
-\pi<-\Omega\leq 
\phi(t)-\phi(s) \leq \Omega<\pi,
\end{equation*}
et, si $t\geq s$, 
\begin{equation*}
\phi(t)-\phi(s)\geq 0,
\end{equation*}
et donc, d'après \eqref{ffvsjhgfss}, $\gamma'(t)\geq 0$.
De même, si $t\leq s$, $\gamma'(t)\leq 0$. Or, on a $\gamma(s)=0$ et donc pour tout $t$, $\gamma(t)\geq 0$. 
$\gamma(t)$ représente la composante du vecteur $X(t)-X(s)$ sur $\sigma(X'(s))$, ce qui nous permet de conclure. 
\ifcase \cras
\end{proof}
\or
\end{preuve}
\fi

\or
La preuve, fondée sur  \eqref{eq20} et \eqref{eq150}, est laissée au lecteur.
\fi


\ifcase \cras
\begin{proof}[Démonstration du lemme \ref{lemmeA}]\
\or
\begin{preuve}[ du lemme   \ref{lemmeA}]\
\fi

\ifcase \cras
\label{lemmeApreuve}
\fi

Soit $X$ 
vérifiant 
\eqref{eq10totnew},
\eqref{eq20},
\eqref{eq100tot},
\eqref{eq110},
\eqref{eq140} (ou \eqref{eq150}).

\begin{figure}[h]    
\psfrag{jt}{$\vec {k}$}
\psfrag{it}{$\vec {\alpha}$}
\psfrag{O}{$O$}
\psfrag{A}{$A$}
\psfrag{B}{$B$}
\psfrag{X}{$X$}
\psfrag{Q}{$Q(s)$}
\psfrag{Xp}{$X'(s)$}
\psfrag{u}{$u$}
\psfrag{v}{$v$}
\begin{center} 
\epsfig{file=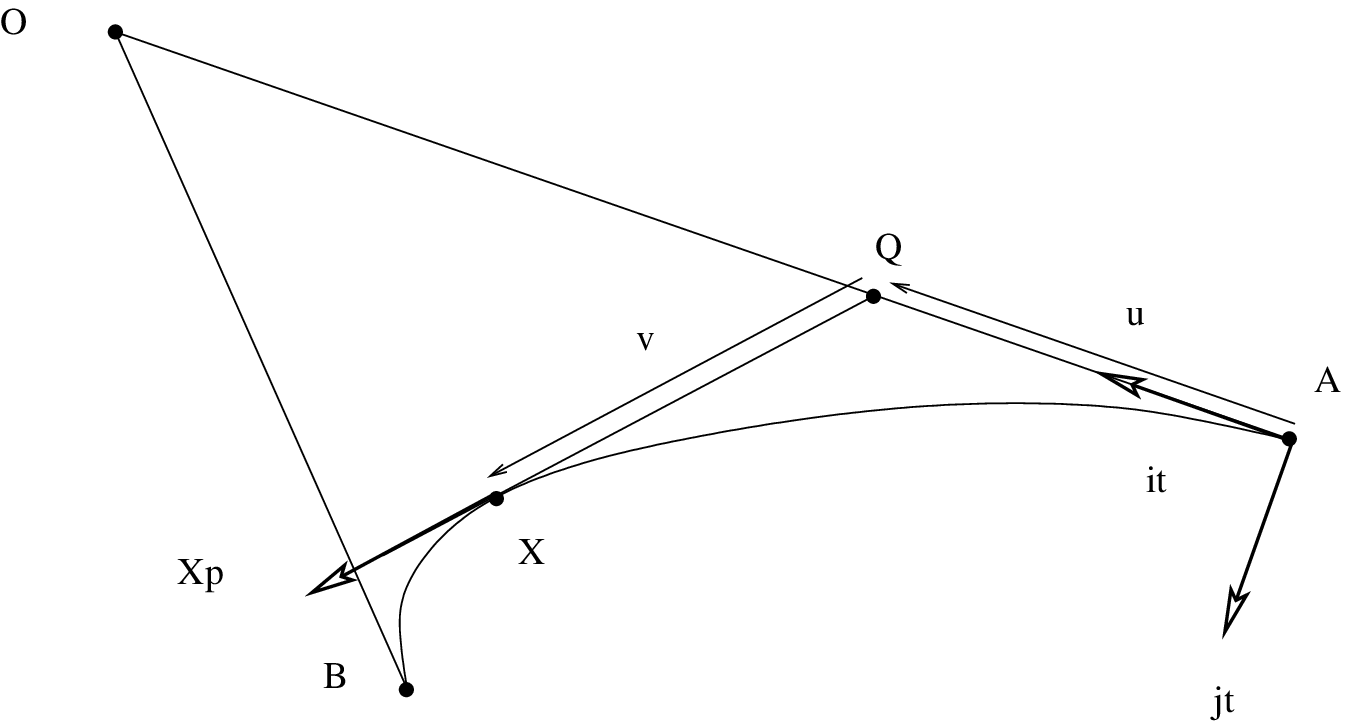, width=9 cm} 
\end{center} 
\caption{\label{longueur_majoreeb}Le repère $\left(A,\vec {\alpha}, \vec {k}\right)$.\ifcase \cras \or /\textit{{The orthonormal frame $\left(A,\vec {i}, \vec {j}\right)$.}}\fi} 
\end{figure}

On considère 
de nouveau le repère orthonormé $\left(A, \vec {\alpha}, \vec {k}\right)$.
On note
de nouveau 
$(  x(s),  y(s))$ les coordonnées de $X$ dans ce repère,
comme indiqué sur la figure \ref{longueur_majoreeb}.
\ifcase \cras
\or
On a les relations habituelles  
\begin{subequations}
\label{eqderphinbbbbnewN}
\begin{align}
\label{eqderphinabnewN}
&\frac{d  x}{ds}=\cos   \phi,\\
\label{eqderphinbbnewN}
&\frac{d y}{ds}=\sin   \phi.
\end{align}
\end{subequations}
\fi
D'après les hypothèses 
\eqref{eq20},
\eqref{eq100d},
\eqref{eq100e} et 
\eqref{eq150}, il existe $s_0\in [0,L[$ tel
que 
\begin{equation}
\label{preuveq00}
\text{$\phi(s)=0$ sur $[0,s]$ et  $0<\phi(s)< \pi$ sur $]s_0,L]$.}
\end{equation}
 Fixons $s\in ]s_0,L]$.
Pour tout point $P=(x_p,y_p)$ de la droite passant par $X(s)$ et porté par $X'(s)$, il existe $\lambda\in \Er$ tel que 
$\overrightarrow{PX(s)}=\lambda X'(s)$, soit d'après 
\ifcase \cras
\eqref{eqderphinbbbb}
\or
\eqref{eqderphinbbbbnewN}
\fi
\begin{equation}
\label{preuveq01}
x(s)-x_p=\lambda  \cos \phi(s),\quad y(s)-y_p=\lambda \sin\phi(s).
\end{equation}
Les deux droites respectives passant par $A$ et portée par $\vec \alpha$
et passant par $X(s)$ et portée par $X'(s)$
se coupent donc un unique point $Q(s)$ d'ordonnée $0$ et d'abscisse donnée par $x_p$
dans \eqref{preuveq01} correspondant à $y_p=0$. On a donc 
\begin{equation}
\label{preuveq10}
\lambda =\frac{y(s)}{\sin\phi(s)},\quad x_p=x(s)-y(s)\frac{ \cos \phi(s)}{\sin\phi(s)}
\end{equation}
Considérons $u(s)$, défini comme  l'abscisse $x_p$ de $Q(s)$ et $v(s)=\lambda$.
Le point $Q(s)$ vérifie donc 
\begin{equation}
\label{preuveq20}
\forall s \in ]s_0,L],\quad
\overrightarrow{AQ(s)}=u(s) \vec \alpha,\quad 
\overrightarrow{Q(s)X(s)}=v(s) \vec X'(s).
\end{equation}
où
\begin{equation}
\label{preuveq30}
u(s)= x(s)-y(s)\frac{ \cos \phi(s)}{\sin\phi(s)}  ,\quad 
v(s)=\frac{y(s)}{\sin\phi(s)} 
\end{equation}
On peut dériver $u$ presque partout et on a, compte tenu  de 
\ifcase \cras
\eqref{eqderphinbbbb}
\or
\eqref{eqderphinbbbbnewN}
\fi
\begin{align}
\text{p.p. sur $]s_0,L[$, }
u'(s)&=x'(s)
-y'(s)\frac{ \cos  \phi(s)}{\sin\phi(s)} 
-y(s)\frac{ -\sin^2\phi(s)\phi'(s)-\cos^2\phi(s)\phi'(s)}{\sin^2\phi(s)} ,\\
&=\cos \phi(s)-\sin\phi(s)\frac{ \cos \phi(s)}{\sin\phi(s)}
+y(s)\frac{ \phi'(s)}{\sin^2\phi(s)} ,\\
\end{align}
et donc 
\begin{equation}
\label{preuveq41}
u'(s)
=y(s)\frac{ \phi'(s)}{\sin^2\phi(s)}.
\end{equation}
Par ailleurs, d'après 
\ifcase \cras
\eqref{eqderphinbbbb},
\or
\eqref{eqderphinbbnewN}
\fi
on a  
\begin{equation*}
y(s)=y(0)+\int_0^s y'(u)du=y(0)+\int_0^s  \sin\phi(u)du.
\end{equation*}
Ainsi, 
d'après 
\eqref{eq20} et 
\eqref{eq160}, $\sin(\phi)\geq 0$ et donc $y(s)\geq 0$, ce qui implique d'après \eqref{eq150} et \eqref{preuveq41} :
\begin{equation}
\label{preuveq50}
\text{$u$ est croissant sur $[s_0,L]$.}
\end{equation}
Ainsi, il existe 
\begin{equation}
\label{preuveq51} 
u(s_0+)=\lim_{\substack{s\to s_0\\ s>s_0}} u(s)\in \{-\infty\}\cap \Er.
\end{equation}
On a 
\begin{equation}
\label{preuveq60}
u(s_0+)\geq 0.
\end{equation}
\begin{figure}[h]    
\psfrag{f1}{$\phi_1$}
\psfrag{jt}{$\vec {k}$}
\psfrag{it}{$\vec {\alpha}$}
\psfrag{O}{$O$}
\psfrag{A}{$A$}
\psfrag{B}{$B$}
\psfrag{X}{$X(s_1)$}
\psfrag{Q}{$Q(s)$}
\psfrag{Xp}{$X'(s_1)$}
\psfrag{u}{$u_1$}
\psfrag{v}{$v_1$}
\begin{center} 
\epsfig{file=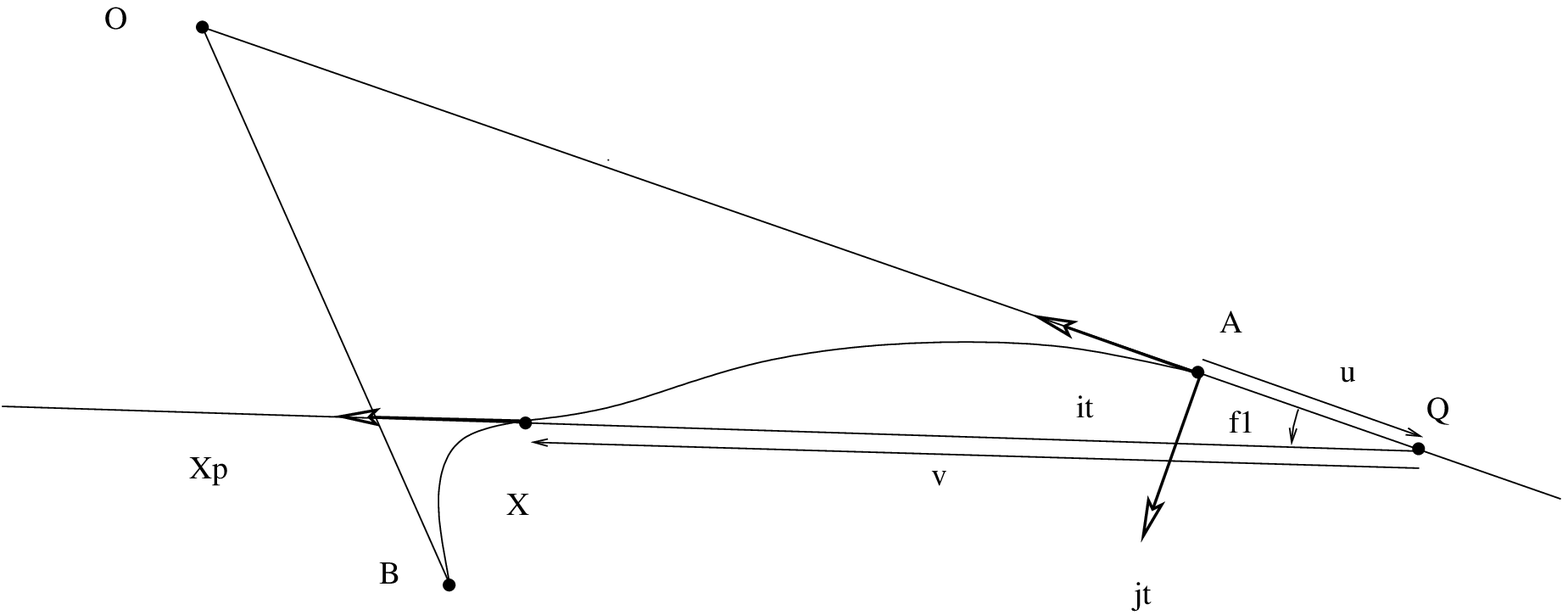, width=12 cm} 
\end{center} 
\caption{\label{longueur_majoree_bis}Le repère $\left(A,\vec {\alpha}, \vec {k}\right)$ avec un $u_1<0$.\ifcase \cras \or /\textit{{The orthonormal frame $\left(A,\vec {i}, \vec {j}\right)$} with $u_1<0$.}\fi} 
\end{figure}
Si ce n'était pas le cas, il existerait $s_1$ tel que 
\begin{equation}
\label{preuveq61}
\phi_1=\phi(s_1)\in ]0,\pi/2[,\quad 
u_1=u(s_1)<0.
\end{equation}
Cela implique que le point $A$ se situe dans le demi-plan ouvert limité par la droite passant par $X(s_1)$, portée par $X'(s_1)$, du côté opposé à la normale 
à la courbe en $X(s_1)$ (voir figure \ref{longueur_majoree_bis}).
Or, d'après le lemme \ref{lemmeconvexite}, la courbe doit aussi se trouver  dans le demi-plan  limité par la droite passant par $X(s_1)$, portée par $X'(s_1)$, du même côté que  la normale 
à la courbe en $X(s_1)$. Ainsi, \eqref{preuveq60} est vrai.
Rappelons que, d'après \eqref{eq100e}, il existe un sous-intervalle $J$ de $[s_0,L]$ tel que 
\begin{equation}
\label{preuveq50b}
\text{$u$ est strictement croissante sur $J$.}
\end{equation}
Ainsi, d'après \eqref{preuveq50}, \eqref{preuveq60} et \eqref{preuveq50b}
\begin{equation}
\label{preuveq70}
u(L)>0.
\end{equation}
Enfin, d'après \eqref{preuveq30}, $v(L)=y(L)/\sin\Omega>0$ et donc 
\begin{equation}
\label{preuveq80}
v(L)>0.
\end{equation}
On conclut en posant $u_0=u(L)$ et $v_0=v(L)$ et en utilisant \eqref{preuveq70} et \eqref{preuveq80} qui impliquent que $O$ est distinct de $A$ et de $B$.
\ifcase \cras
\end{proof}
\or
\end{preuve}
\fi


\newcommand{\preuvelocaleA}{%

\begin{enumerate}

\item
Il suffit de choisir $\mu=AB$, non nul.

\item
Montrons maintenant l'existence de $\nu$.

\begin{figure}[h]
\psfrag{jt}{$\vec {j}$}
\psfrag{it}{$\vec {i}$}
\psfrag{O}{$O$}
\psfrag{A}{$A$}
\psfrag{B}{$B$}
\psfrag{a}{$\alpha$}
\psfrag{b}{$\beta$}
\psfrag{X}{$X$}
\psfrag{xt}{$\widetilde x$}
\psfrag{bt}{$\widetilde b$}
\begin{center} 
\epsfig{file=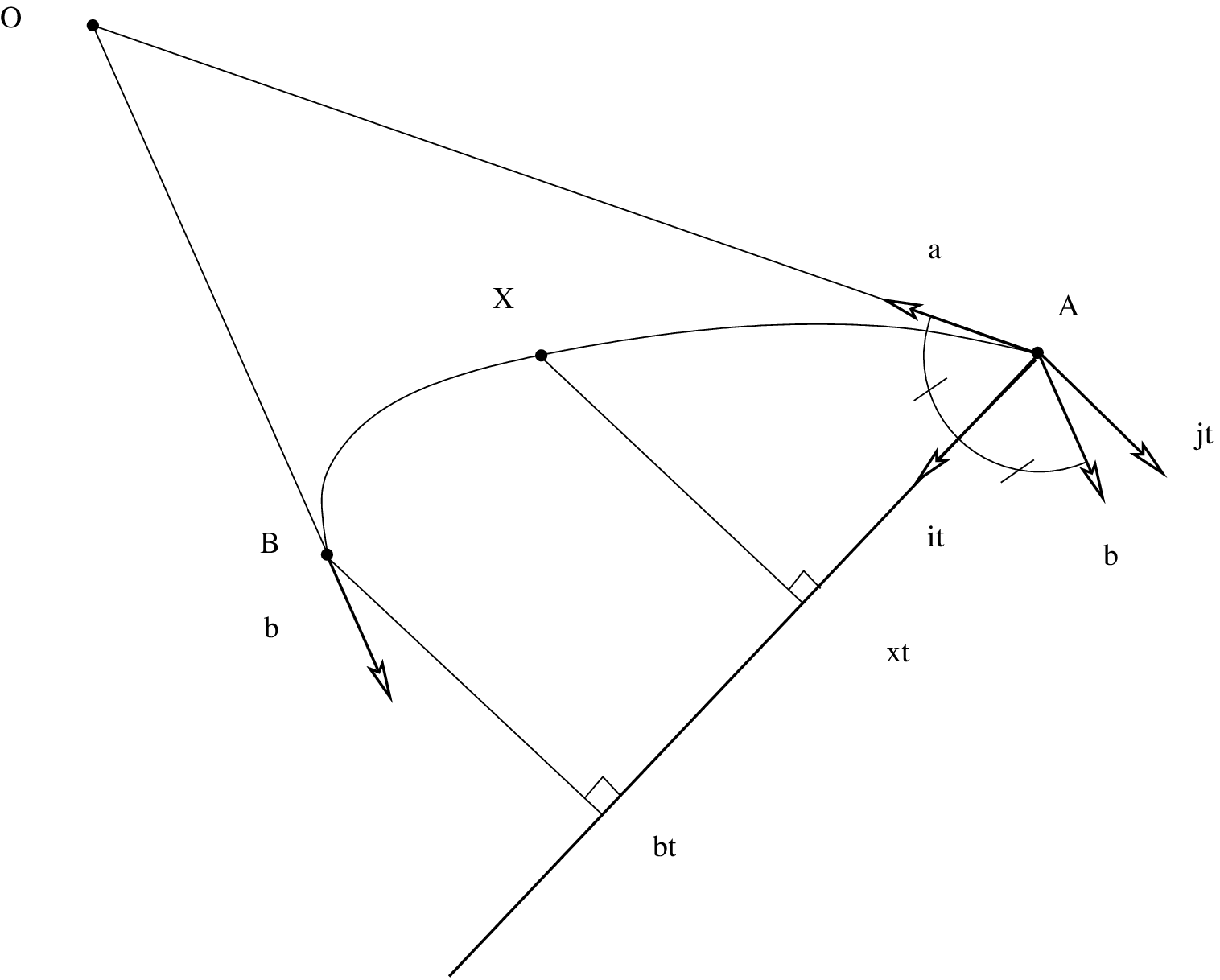, width=9 cm} 
\end{center} 
\caption{\label{longueur_majoree}Le repère $\left(A,\vec {i}, \vec {j}\right)$.\ifcase \cras \or /\textit{{The orthonormal frame $\left(A,\vec {i}, \vec {j}\right)$}}\fi}
\end{figure}
Comme indiqué sur la figure \ref{longueur_majoree}, on considère le repère orthonormé $\left(A, \vec {i}, \vec {j}\right)$
défini par 
\begin{equation}
\label{eq270}
\left( \widehat{{\vec \alpha},\vec {i}} \right)=\frac{\Omega}{2}.
\end{equation}
et on définit  l'angle $\widetilde \phi$ analogue à celui défini par \eqref{eq120}
\begin{equation}
\label{eq30b}
\widetilde \phi=
\left( \widehat{\vec i,X'(s)}\right).
\end{equation}
Ainsi, dans ce repère, en notant $X(s)=(\widetilde x(s),\widetilde y(s))$, les coordonnées des points de  la courbe $X$, 
on a les relations habituelles 
identiques à 
\ifcase \cras
\eqref{eqderphinbbbb}
\or
\eqref{eqderphinbbbbnewN}
\fi
 :
\begin{subequations}
\label{eqderphin}
\begin{align}
\label{eqderphina}
&\frac{d\widetilde x}{ds}=\cos  \widetilde\phi,\\
\label{eqderphinb}
&\frac{d\widetilde y}{ds}=\sin  \widetilde\phi.
\end{align}
\end{subequations}
On a donc, d'après  \eqref{eq270} et \eqref{eq30b}
\begin{equation*}
\widetilde \phi=
\left( \widehat{\vec i,X'(s)}\right)=
\left( \widehat{\vec i,\vec \alpha}\right)+
\left( \widehat{\vec \alpha,X'(s)}\right),
\end{equation*}
et donc 
\begin{equation}
\label{eq280}
\widetilde \phi=
\phi-\frac{\Omega}{2}.
\end{equation}
L'hypothèse \eqref{eq160} implique donc que
\begin{equation}
\label{eq299}
\widetilde \phi\in \left[-\frac{\Omega}{2},\frac{\Omega}{2}\right].
\end{equation}
D'après \eqref{eq20}, on a donc $\Omega/2\in [0,\pi/2[$ et donc  
\begin{equation*}
\cos \widetilde \phi \geq \cos(\Omega/2)>0,
\end{equation*}
et donc 
\begin{equation}
\label{eq290}
0< \frac{1}{\cos \widetilde \phi}\leq \frac{1}{\cos(\Omega/2)}.
\end{equation}
Enfin, d'après \eqref{eqderphina} et \eqref{eq290}, on a pour toute courbe $X\in \mathcal{E}$, en notant $\widetilde b$ l'abscisse de $B$ dans le repère 
$\left(A, \vec {i}, \vec {j}\right)$ 
\begin{equation}
\label{eq291}
L=\int_0 ^L ds=
\int_0 ^{\widetilde b} 
\frac{d\widetilde x}{\cos\widetilde \phi},
\end{equation}
et donc 
\begin{equation*}
L
\leq  
\frac{1}{\cos(\Omega/2)}
\int_0 ^{\widetilde b} 
d\widetilde x
\end{equation*}
et donc 
\begin{equation}
\label{eq300}
L\leq 
\frac{\widetilde b}{\cos(\Omega/2)},
\end{equation}
qui ne dépend que des trois points $O$, $A$ et $B$.
Notons que $\mathcal{E}$ est non vide ; en effet, 
il contient par exemple une parabole (voir la parabole introduite en section \ref{introduction}, qui est de courbure algébrique de signe constant). 
On peut aussi utiliser la courbe construite en section \ref{uniquecercledroite}.
Ainsi,  il existe une courbe avec $L$ non nul, d'après la première inégalité de \eqref{lemmeBeq01},  et    la constante $\frac{\widetilde b}{\cos(\Omega/2)}$  est strictement positive. 
\end{enumerate}
}

\newcommand{\preuvelocaleB}{%
\begin{equation*}
\Omega=\left|\Omega\right|=\left|\phi(B)-\phi(A)\right|=\left|\int_A ^B d\phi\right|=\left|\int_0 ^L   cds\right|\leq \int_0 ^L |c| ds\\
\leq \vnorm[L^{\infty}(0,L)]{c} \int_0 ^L   ds = \vnorm[L^{\infty}(0,L)]{c} L,
\end{equation*}
et donc 
\begin{equation}
\label{fjfhfhf}
\Omega\leq\ \vnorm[L^{\infty}(0,L)]{\vnorm{X''}} L.
\end{equation}
On conclut ce point 
grâce à  \eqref{lemmeBeq01} 
 et en posant $\delta=\Omega/\nu$.
}


\ifcase \faussepreuve

\ifcase \cras

\begin{proof}[Démonstration de la proposition \ref{existencetheoprop01new}]\

\ifcase \cras
\label{preuveexistencetheoprop01new}
\fi

\begin{enumerate}

\item

\preuvelocaleA

\item
 On écrit
\preuvelocaleB
\end{enumerate}

\end{proof}


\begin{proof}[Démonstration du théorème \ref{existencetheoprop01juste}]\

\ifcase \cras
\label{existencetheoprop01justepreuve}
\fi

On  a déjà remarqué, dans la preuve de la proposition \ref{existencetheoprop01new}
que $\mathcal{E}$ était non vide.

D'après \eqref{existencetheolem02dfg}, on peut poser 
\begin{equation} 
\label{existencetheolem02dfgbis}
\delta_0=\inf_{Z \in \mathcal{E}} \vnorm[L^{\infty}(0,L(Z))]{\vnorm{Z''}}>0.
\end{equation}
Pour toute la suite, on peut considérer  la fonction $g$ de $\mathcal{E}$ dans $\Erps$ définie par 
\begin{equation}
\label{pdueq01}
\forall X\in \mathcal{E},\quad
g(X)=\vnorm[L^{\infty}(0,L(X))]{\vnorm{X''}}.
\end{equation}
D'après la  définition \eqref{existencetheolem02dfgbis} de $\delta_0>0$,
il existe une suite $(X_n)$ de $\mathcal{E}$ telle que 
\begin{equation}
\label{discretcaspthenconeq70gdgd}
\lim_{n\to \infty} g(X_n)=\delta_{0}.
\end{equation}
Montrons que $X_n$ admet une suite extraite convergeant vers un élément $X$ de $\mathcal{E}$,
ce qui nous permettra de conclure.
Il nous faut tout d'abord
construire une sous-suite de $X_n$ telle que la suite des longueurs $L(X_n)$ converge vers
un réel $L_0 > 0$.
D'après \eqref{lemmeBeq01}, on a 
\begin{equation}
\label{existencetheoprop01justepreuveeq00}
\forall n\in \En,\quad
\mu \leq L(X_n) \leq \nu
\end{equation}
Ainsi,  puisque $L(X_n)$ appartient à un compact de $\Er$, 
on peut extraire de $X_n$ une sous suite, encore notée de la même façon, de telle sorte que la suite 
$L(X_n)$ converge  vers un réel $L_0>0$ : 
\begin{equation}
\label{discretcaspthenconeq70gdgdbis}
\lim_{n\to \infty} L(X_n)=L_0.
\end{equation}
On est donc de nouveau dans le cadre de la preuve de \cite[Proposition 1]{MR0089457},
que l'on adapte à notre cas.
En effet, la preuve de \cite[Proposition 1]{MR0089457} est fondée sur le fait que $L_0$ est la borne inférieure des longueurs $L(X_n)$, puisque Dubins
cherche  une courbe minimisant la longueur. Ici ce n'est plus le cas et on nous allons d'abord nous ramener à un intervalle de référence de longueur $1$ 
avant d'extraire une sous-suite de $X_n$.
On considère donc la suite de fonctions $\widetilde X_n$ définies par 
\begin{equation}
\label{existencetheoprop01justepreuveeq01}
\forall t\in [0,1],\quad
\widetilde X_n(t)=X_n\left(L(X_n) t\right).
\end{equation}
Nous allons alors extraire de $\widetilde X_n$ une sous-suite convergente.
Il est clair que, puisque $X_n$ appartient à $\mathcal{E}$,  on a, pour tout $n$ : 
\begin{equation}
\label{eq100totexistencetheoprop01justepreuveeq20}
\widetilde X_n\in W^{2,\infty}(0,1;\Er^2),
\end{equation}
avec 
\begin{subequations}
\label{eq100totexistencetheoprop01justepreuveeq30}
\begin{align}
\label{eq100totexistencetheoprop01justepreuveeq30a}
&\forall t\in [0,1],\quad \widetilde X'_n(t)=L(X_n) X_n'\left(L(X_n) t\right),\\
\label{eq100totexistencetheoprop01justepreuveeq30b}
&\text{p.p. sur } ]0,1[, \quad \widetilde X ''_n(t)=L^2(X_n) X_n''\left(L(X_n) t\right).
\end{align}
\end{subequations}
et donc que 
\begin{subequations}
\label{eq100totexistencetheoprop01justepreuveeq10}
\begin{align}
\label{eq100totexistencetheoprop01justepreuveeq10a}
&\forall t\in [0,1],\quad \vnorm{X'_n(s)}=L(X_n),\\
\label{eq100totexistencetheoprop01justepreuveeq10b}
&\widetilde X_n(0)=A,\\
\label{eq100totexistencetheoprop01justepreuveeq10c}
&\widetilde X_n(1)=B,\\
\label{eq100totexistencetheoprop01justepreuveeq10d}
&\widetilde X'_n(0)=L(X_n)\vec \alpha,\\
\label{eq100totexistencetheoprop01justepreuveeq10e}
&\widetilde X'_n(1)=L(X_n) \vec \beta.
\end{align}
\end{subequations}
On considère une détermination continue de l'angle 
$ \widetilde \phi_n$ défini par 
\begin{equation}
\label{eq100totexistencetheoprop01justepreuveeq50}
\forall t\in [0,1],\quad 
  \widetilde \phi_n(t)=
\left( \widehat{\vec \alpha,\widetilde X'_n(t)} \right), 
\end{equation}
et on a, puisque l'angle $\phi_n$ associé à $X_n$ est croissant,  
\begin{equation}
\label{eq100totexistencetheoprop01justepreuveeq60}
\text{$\widetilde \phi_n$ est croissant.}
\end{equation}
Remarquons que 
\eqref{existencetheoprop01justepreuveeq00}
et \eqref{eq100totexistencetheoprop01justepreuveeq10a}
impliquent 
\begin{equation}
\label{eq100totexistencetheoprop01justepreuveeq100}
\forall n\in \En,\quad 
\forall t\in [0,1],\quad \vnorm{\widetilde  X'_n(t)} \leq \nu.
\end{equation}
et donc 
la suite $X'_n$ est une famille de fonctions uniformément bornées. 
De plus, on a 
\begin{equation*}
\forall n\in \En, \quad 
\forall t_1,t_2\in [0,1],\quad 
\vnorm{\widetilde X'_n(t_1)-\widetilde X'_n(t_2)}
\leq \vnorm[L^{\infty}(0,1)]{\vnorm{\widetilde X''_n}}|t_1-t_2|,
\end{equation*}
soit d'après \eqref{eq100totexistencetheoprop01justepreuveeq30b} : 
\begin{equation}
\label{existencetheoprop01justepreuveeq110}
\forall n\in \En, \quad  
\forall t_1,t_2\in [0,1],\quad 
\vnorm{\widetilde X'_n(t_1)-\widetilde X'_n(t_2)}
\leq L^2(X_n) g(X_n)|t_1-t_2|,
\end{equation}
et donc, en posant $K=\nu ^2 \sup_{n} g(X_n)<+\infty$, car la suite $g(X_n)$ converge et en utilisant \eqref{existencetheoprop01justepreuveeq00}, on a 
\begin{equation}
\label{pdueq21}
\forall n,
\quad
\forall t_1,t_2\in [0,1],\quad
\vnorm{\widetilde X'_n(t_1)-X\widetilde '_n(t_2)}
\leq 
K|t_1-t_2|.
\end{equation}
Ainsi, la famille $X'_n$ est une famille uniformément bornée et équicontinue   et,
d'après le théorème d'Ascoli%
\footnote{voir par exemple \url{http://www.bibmath.net/dico/index.php?action=affiche&quoi=./e/equicontinu.html}}, 
il existe une suite extraite de $\widetilde X_n$, encore notée de la même façon, dont la dérivée converge uniformément sur $[0,1]$ vers une fonction $\widetilde Y$ continue. 
Puisque 
\begin{equation*}
\widetilde X_n(t)=\widetilde X_n(0)+\int_0 ^t \widetilde X'_n(u)du=
A+\int_0 ^t \widetilde X'_n(u)du,
\end{equation*}
on a, pour tout $t\in [0,1]$,
\begin{equation*}
\vnorm{\widetilde  X_n(t)-\widetilde  X_m(t)}\leq 
\int_0 ^t \vnorm{\widetilde  X'_n(u)-\widetilde X'_m(u)}du
\leq
\sup_{u\in [0,1]}  \vnorm{\widetilde  X'_n(u)-\widetilde  X'_m(u)},
\end{equation*}
et donc la suite $\widetilde X_n$ converge uniformément sur $[0,1]$ vers une fonction $\widetilde X$. Puisque $\widetilde X'_n$ converge uniformément vers 
$\widetilde Y$ et que  $\widetilde X_n$ converge vers $\widetilde X$, on 
a ${\left(\widetilde X\right)}'={\left(\lim \widetilde X_n\right)}'=\lim \widetilde X'_n=\widetilde Y$.  Ainsi, $\widetilde X$ est de classe ${\mathcal{C}}^1$.
Par ailleurs,  de \eqref{pdueq21}, on déduit, en passant à la limite $n$ tendant vers l'infini : 
\begin{equation}
\label{pdueq22}
\forall t_1,t_2\in [01],\quad
\vnorm{\widetilde X'(t_1)-\widetilde X'(t_2)}
\leq 
K|t_1-t_2|,
\end{equation}
et donc, d'après par exemple,  
\cite[Proposition VIII.3 et Corollaire VIII.4]{MR697382},
$\widetilde X$ appartient à 
$W^{2,\infty}(0,1;\linebreak[1] \Er^2)$.
Si on passe à la limite $n$ tendant vers l'infini dans \eqref{existencetheoprop01justepreuveeq110}, on a aussi 
\begin{equation}
\label{existencetheoprop01justepreuveeq120}
\forall t_1,t_2\in [0,1],\quad 
\vnorm{\widetilde X'(t_1)-\widetilde X'(t_2)}
\leq L_0^2 \delta_0|t_1-t_2|,
\end{equation}
et  donc de nouveau en 
utilisant{\footnote{voir aussi \texttt{https://fr.wikipedia.org/wiki/Application\udsc lipschitzienne\#cite\udsc note-Wikiversité}}}
\cite[Proposition VIII.3 et Corollaire VIII.4]{MR697382} et \cite[Théorème 8.18]{MR662565}, on a 
\begin{equation}
\label{existencetheoprop01justepreuveeq130}
\vnorm[L^{\infty}(0,1)]{\vnorm{\widetilde X''}}\leq L_0^2 \delta_0.
\end{equation}
Revenons enfin à l'intervalle $[0,L_0]$ en considérant la fonction $X \in W^{2,\infty}(0,L_0;\linebreak[1] \Er^2)$ définie par 
\begin{equation}
\label{existencetheoprop01justepreuveeq01ret}
\forall s\in [0,L_0],\quad
X(s)=\widetilde X\left(\frac{s}{L_0}\right).
\end{equation}

Montrons, pour conclure, que $X$ appartient à $\mathcal{E}$ et minimise $g(Z)$ pour $Z$ décrivant $\mathcal{E}$.
Il est clair qu'on a 
\begin{equation*}
\forall s\in [0,L_0],\quad
X'(s)=\frac{1}{L_0}\widetilde X'\left(\frac{s}{L_0}\right).
\end{equation*}
et qu'en passant à la limite $n\to\infty$  dans \eqref{eq100totexistencetheoprop01justepreuveeq30a} et 
\eqref{eq100totexistencetheoprop01justepreuveeq10}, on obtient donc 
\begin{align*}
&\forall s \in [0,L_0], \quad \vnorm{X'(s)}=1,\\
&  X(0)=A,\\
&  X(L_0)=B,\\
&  X'(0)= \vec \alpha,\\
&  X'(L_0)=  \vec \beta.
\end{align*}
En considérant  une détermination continue de l'angle 
$\phi$ défini par 
\begin{equation*}
\forall s\in [0,L_0],\quad 
 \phi(s)=
\left( \widehat{\vec \alpha, X'(s)} \right)=
\lim_{n\to \infty}\left( \widehat{\vec \alpha, \widetilde X'_n\left(\frac{s}{L_0}\right)} \right)=
\lim_{n\to \infty}   \widetilde \phi_n  \left( \frac{s}{L(X_n)}   \right).
\end{equation*}
et on a, selon \eqref{eq100totexistencetheoprop01justepreuveeq60}, 
\begin{equation*}
\text{$\phi$ est croissant.}
\end{equation*}
Naturellement, on vérifie que 
\begin{equation*}
L(X)=\int_0 ^{L_0}\vnorm{X'(s)}ds,
\end{equation*}
soit 
\begin{equation}
\label{existencetheoprop01justepreuveeq250}
L(X)=L_0.
\end{equation}
Ainsi
 $  X$ appartient à 
$W^{2,\infty}(0,L_0;\Er^2)$ et à $\mathcal{E}$.
Enfin, 
on a, d'après \eqref{existencetheoprop01justepreuveeq130}, 
\begin{equation*}
\vnorm[L^{\infty}(0,L_0)]{\vnorm{X''}}=\frac{1}{L_0^2}\vnorm[L^{\infty}(0,1)]{\vnorm{\widetilde X''}}
\leq \frac{1}{L_0^2} L_0^2 \delta_0=\delta_0,
\end{equation*}
et donc 
\begin{equation}
\label{existencetheoprop01justepreuveeq150}
\vnorm[L^{\infty}(0,L(X))]{\vnorm{X''}} \leq \delta_0.
\end{equation}
D'après la définition  \eqref{existencetheolem02dfgbis} de $\delta_0$, puisque $X$ appartient à $\mathcal{E}$, on a donc 
\begin{equation}
\label{existencetheoprop01justepreuveeq160}
\vnorm[L^{\infty}(0,L(X))]{\vnorm{X''}} \geq \delta_0.
\end{equation}
Finalement, grâce à \eqref{existencetheoprop01justepreuveeq150} et \eqref{existencetheoprop01justepreuveeq160}, on a 
\begin{equation*}
\vnorm[L^{\infty}(0,L(X))]{\vnorm{X''}} =\delta_0.
\end{equation*}
Ainsi, $X$ appartient à $\mathcal{E}$ et $X$ est de maximum de courbure minimal, ce qui permet de conclure cette preuve.
\end{proof}

\fi

\or


\ifcase \cras
\begin{proof}[Démonstration du théorème  \ref{existencetheoprop01}]\
\or
\begin{preuve}[ du théorème  \ref{existencetheoprop01}]\
\fi

\textbf{Attention, ce résultat est faux !!!! Néanmoins, les deux estimations \eqref{lemmeBeq01} (seule la majoration est conservée) et \eqref{existencetheolem02dfg} sont correctes et seront utilisées dans le cas discret dans \cite{piece6_optimale_JB_2019_X4}. Et puis résultat rétabli !!!!!}

\label{refpreuvefausse}

Nous procédons en plusieurs étapes. 

\begin{enumerate}

\item
L'ensemble $\mathcal{E}$ est non vide puisqu'il contient par exemple une parabole (voir la parabole introduite en section \ref{introduction}, qui est de courbure algébrique de signe constant). 
On peut aussi utiliser la courbe construite en section \ref{uniquecercledroite}.

\item
Montrons maintenant qu'il  existe deux constantes strictement positives $\mu$ et $\nu$, ne dépendant que des points $O$, $A$ et $B$ vérifiant \eqref{lemmeBeq01}. 

Il  suffit de choisir $\mu=AB$.

\preuvelocaleA

\ifcase \faussepreuve
On peut vérifier que $\mathcal{E}$ est non vide puisqu'il contient par exemple une parabole (voir la parabole introduite en section \ref{introduction}, qui est de courbure algébrique de signe constant). 
On peut aussi utiliser la courbe construite en section \ref{uniquecercledroite}. Ainsi, 
\or
Puisque $\mathcal{E}$ est non vide, 
\fi
il existe une courbe avec $L$ non nul  et    la constante $\frac{\widetilde b}{\cos(\Omega/2)}$  est strictement positive.

\item
Ensuite, on montre qu'il  existe $\delta>0$ telle que pour toute courbe $X$ de
 $\mathcal{E}$, on ait  \eqref{existencetheolem02dfg}. 
Pour cela, on écrit
\preuvelocaleB

\item
\ifcase \cras
Pour toute la suite, on considère la fonction $g$ de $\mathcal{E}$ dans $\Erps$ définie par 
\begin{equation}
\label{pdueq01}
\forall X\in \mathcal{E},\quad
g(X)=\vnorm[L^{\infty}(0,L(X))]{\vnorm{X''}}.
\end{equation}
\or
On rappelle la définition \eqref{pdueq01} de $g$.
\fi
D'après la  définition \eqref{existencetheolem02dfgbis} de $\delta_0>0$,
il existe une suite $(X_n)$ de $\mathcal{E}$ telle que 
\begin{equation}
\label{discretcaspthenconeq70gdgd}
\lim_{n\to \infty} g(X_n)=\delta_{0}.
\end{equation}
Montrons que $X_n$ admet une suite extraite convergeant vers un élément $X$ de $\mathcal{E}$,
ce qui nous permettra de conclure.
Il nous faut tout d'abord
construire une sous-suite de $X_n$ telle que les longueurs $L(X_n)$ décroissent vers
un réel $L_0 > 0$.
D'après \eqref{lemmeBeq01}, on a 
\begin{equation}
\forall n\in \En,\quad
L(X_n)\geq \mu.
\end{equation}

\textbf{Attention, ici se trouve la faute, vue par le rapporteur !}. Il écrit dans son rapport : "p11 l22-23 Je ne comprends pas l'argument : le fait que la suite des $L(Xn)$
soit minorée par $m$ n'implique pas que l'on puisse extraire une sous-suite 
$ L(X_{\phi(n)})$ qui soit décroissante vers un $L_0\geq m$."

\textit{En effet, j'ai été trop vite et ai voulu absolument calquer la preuve de Dubins. On peut seulement affirmer qu'il existe un élement de l'ensemble des $X_n$ qui converge vers $\delta_0$, mais cet élément 
ne permet pas d'extraire une sous-suite !!!! }

\textit{On peut mofier le raisonnement utiliser l'estimation \eqref{lemmeBeq01} pour montrer, les valeurs de $L(X_n)$ étant dans un compact, extraire une sous-suite qui converge vers un certain $L_0$, qui peut
être choisi comme la plus petite des valeurs d'adhérences, mais dans ce cas $L_0$ n'est plus nécessairement un minorant des $L(X_n)$ et dans ce cas, le raisonnement de Dubins ne peut plus s'appliquer.}

\textit{Conclusion : 
\begin{itemize}
\item
Ouf, je ne publie pas de conneries, finalement !
\item
On fait un travail différent de Dubins !
\ifcase \cras
\or
\item
Ce  théorème  \ref{existencetheoprop01} n'est pas finalement utile, grâce au théorème \ref{existenceunicitetheoprop01}, lui, nouveau !
\fi
\item 
deux estimations resrent vraies, présentées dans \cite{piece6_optimale_JB_2019_X4}.
\item
C'est légtime en fait ; je cherchais à la fois une courbe minimisant le maximum de la courbure et minisant la longueur ! On ne peut pas tout avoir, même si finalement on a tout !
Car notre problème  a la même solution que Dubins avec $R=R_a(O,A,B)$.
\item
Vient d'être rétabli !!!
\end{itemize}}

Ainsi,  
on peut extraire de $X_n$ une sous suite, encore notée de la même façon, de telle sorte que la suite 
$L(X_n)$ converge en décroissant vers un réel $L_0>0$. 
On est donc de nouveau dans le cadre de la preuve de \cite[Proposition 1]{MR0089457},
que l'on adapte à notre cas.
Puisque $\vnorm{X'_n}=1$, la suite $X'_n$ est une famille de fonctions uniformément bornées. De plus, on a pour tout $n$, 
pour tout $s_1,s_2$ dans $[0,L(X_n)]$,
\begin{equation*}
\vnorm{X'_n(s_1)-X'_n(s_2)}
\leq \vnorm[L^{\infty}(0,L(X_n))]{\vnorm{X''_n}}|s_1-s_2|
=g(X_n)|s_1-s_2|,
\end{equation*}
et donc, en posant $K=\sup_{n} g(X_n)<+\infty$, car la suite $g(X_n)$ converge, on a 
\begin{equation}
\label{pdueq20}
\forall n,
\quad
\forall s_1,s_2\in [0,L(X_n)],\quad
\vnorm{X'_n(s_1)-X'_n(s_2)}
\leq 
K|s_1-s_2|,
\end{equation}
Puisque $[0,L_0]\subset [0,L(X_n)]$, on a \textit{a fortiori}
\begin{equation}
\label{pdueq21}
\forall n,
\quad
\forall s_1,s_2\in [0,L_0],\quad
\vnorm{X'_n(s_1)-X'_n(s_2)}
\leq 
K|s_1-s_2|,
\end{equation}
Ainsi, la famille $X'_n$ est une famille équicontinue et,
d'après le théorème d'Ascoli, 
il existe une suite extraite de $X_n$, encore notée de la même façon, dont la dérivée converge uniformément sur $[0,L_0]$ vers une fonction $Y$ de classe ${\mathcal{C}}^1$. 
Puisque 
\begin{equation*}
X_n(s)=X_n(0)+\int_0 ^ s X'_n(t)dt=
A+\int_0 ^s X'_n(t)dt,
\end{equation*}
on a, pour tout $s\in [0,L_0]$,
\begin{equation*}
\vnorm{X_n(s)-X_m(s)}\leq 
\int_0 ^s \vnorm{X'_n(t)-X'_m(t)}dt
\leq
L_0 \sup_{t\in [0,L_0]}  \vnorm{X'_n(t)-X'_m(t)},
\end{equation*}
et donc la suite $X_n$ converge uniformément sur $[0,L_0]$ vers une fonction $X$. Puisque $X'_n$ converge uniformément vers 
$Y$ et que  $X_n$ converge vers $X$, on 
a $X'=(\lim X_n)'=\lim X'_n=Y$.  Ainsi, $X$ est de classe ${\mathcal{C}}^1$.
Ainsi, de \eqref{pdueq21}, on déduit, en passant à la limite $n$ tendant vers l'infini : 
\begin{equation}
\label{pdueq22}
\forall s_1,s_2\in [0,L_0],\quad
\vnorm{X'(s_1)-X'(s_2)}
\leq 
K|s_1-s_2|,
\end{equation}
et donc, d'après 
\cite[Proposition VIII.3 et Corollaire VIII.4]{MR697382},
$X$ appartient à 
$X\in W^{2,\infty}(0,L;\Er^2)$.
On a aussi 
\begin{equation*}
X(0)=\lim_{n \to \infty} X_n(0)=\lim_{n \to \infty} A=A,\quad 
X'(0)=\lim_{n \to \infty} X'_n(0)=\lim_{n \to \infty} \vec \alpha=\vec \alpha.
\end{equation*}
Il est élémentaire de continuer la preuve en montrant que $X(L_0)=B$ et $X'(L_0)=\vec \beta$, puisque
\begin{multline*}
\vnorm{X(L_0)-B}
\leq \vnorm{X(L_0)-X_n(L_0)}+\vnorm{X_n(L_0)-B}
    =\vnorm{X(L_0)-X_n(L_0)}+\vnorm{X_n(L_0)-X_n(L(X_n))},\\
\leq\vnorm{X(L_0)-X_n(L_0)}+\int_{L_0}^{L(X_n)}\vnorm{X_n'(s)}ds
\leq\vnorm{X(L_0)-X_n(L_0)}+L(X_n)-L_0,
\end{multline*}
qui tend vers zéro quand $n$ tend vers l'infini.
On a aussi, d'après \eqref{pdueq20}
\begin{multline*}
\vnorm{X'(L_0)-\vec \beta}
\leq \vnorm{X'(L_0)-X_n'(L_0)}+\vnorm{X_n'(L_0)-\vec \beta},\\
=\vnorm{X'(L_0)-X'_n(L_0)}+\vnorm{X'_n(L_0)-X'_n(L(X_n))}
\leq
\vnorm{X'(L_0)-X'_n(L_0)}+K|L_0-L(X_n)|,
\end{multline*}
qui tend vers zéro quand $n$ tend vers l'infini.
Il faut compléter la preuve de Dubins en montrant 
que l'hypothèse \eqref{eq140} (ou \eqref{eq150}) est vérifiée. 
On a
\begin{equation*}
 \forall s\in [0,L_0],\quad
 \phi(s)=\left( \widehat{\vec \alpha,X'(s)} \right)
 =\lim_{n\to \infty} \left( \widehat{\vec \alpha,X_n'(s)} \right)=\lim_{n\to \infty} \phi_n(s).
 \end{equation*}
Puisque chacune des fonctions $\phi_n$ est croissante, la limite l'est encore et donc, à la limite,  
\eqref{eq150} est vérifiée.
Ainsi, $X$ appartient à $\mathcal{E}$ et $X$ est de maximum de courbure minimal. 
\end{enumerate}
\ifcase \cras
\end{proof}
\or
\end{preuve}
\fi
\fi


\ifcase \cras

\ifcase \cras
\begin{remark}
\or
\begin{remarque}
\fi
\label{fdddsfsdf}
Les preuves  de la proposition 
\ifcase \faussepreuve
\ref{existencetheoprop01new}
\or
\ref{existencetheoprop01}
\fi
et du théorème 
\ifcase \faussepreuve
\ref{existencetheoprop01juste} 
\or
\textbf{faux} \ref{existencetheoprop01}
\fi
sont fondées sur les deux inégalités
\eqref{eq299} et
\eqref{eq290}
qui proviennent de \eqref{eq160}, elle-même conséquence de 
\eqref{eq150}. 
L'autre inégalité fondamentale de ces preuves
\eqref{fjfhfhf} est encore valable sans l'hypothèse \eqref{eq150}.
Autrement dit, la
proposition  
\ifcase \faussepreuve
\ref{existencetheoprop01new}
\or
\ref{existencetheoprop01}
\fi
et le théorème 
\ifcase \faussepreuve
\ref{existencetheoprop01juste} 
\or
\textbf{faux} \ref{existencetheoprop01}
\fi
sont encore valables si l'on remplace l'hypothèse \eqref{eq150}
par l'hypothèse plus générale \eqref{eq160}.
\ifcase \cras
\end{remark}
\or
\end{remarque}
\fi

\fi


\ifcase \cras
\begin{proof}[Démonstration  du lemme \ref{uniquecercledroitelem01}]\
\or
\begin{preuve}[  du lemme \ref{uniquecercledroitelem01}]\
\fi

\ifcase \cras
\label{uniquecercledroitelem01preuve}
\fi

\begin{figure}[h]    
\psfrag{A}{$A$}
\psfrag{B}{$B$}
\psfrag{O}{$O$}
\psfrag{a}{$\alpha$}
\psfrag{b}{$\beta$}
\psfrag{d1}{$d_1$}
\psfrag{d2}{$d_2$}
\psfrag{C1}{$C_1$}
\psfrag{C2}{$C_2$}
\psfrag{c1}{$\mathcal{C}_1$}
\psfrag{c2}{$\mathcal{C}_2$}
\begin{center} 
\epsfig{file=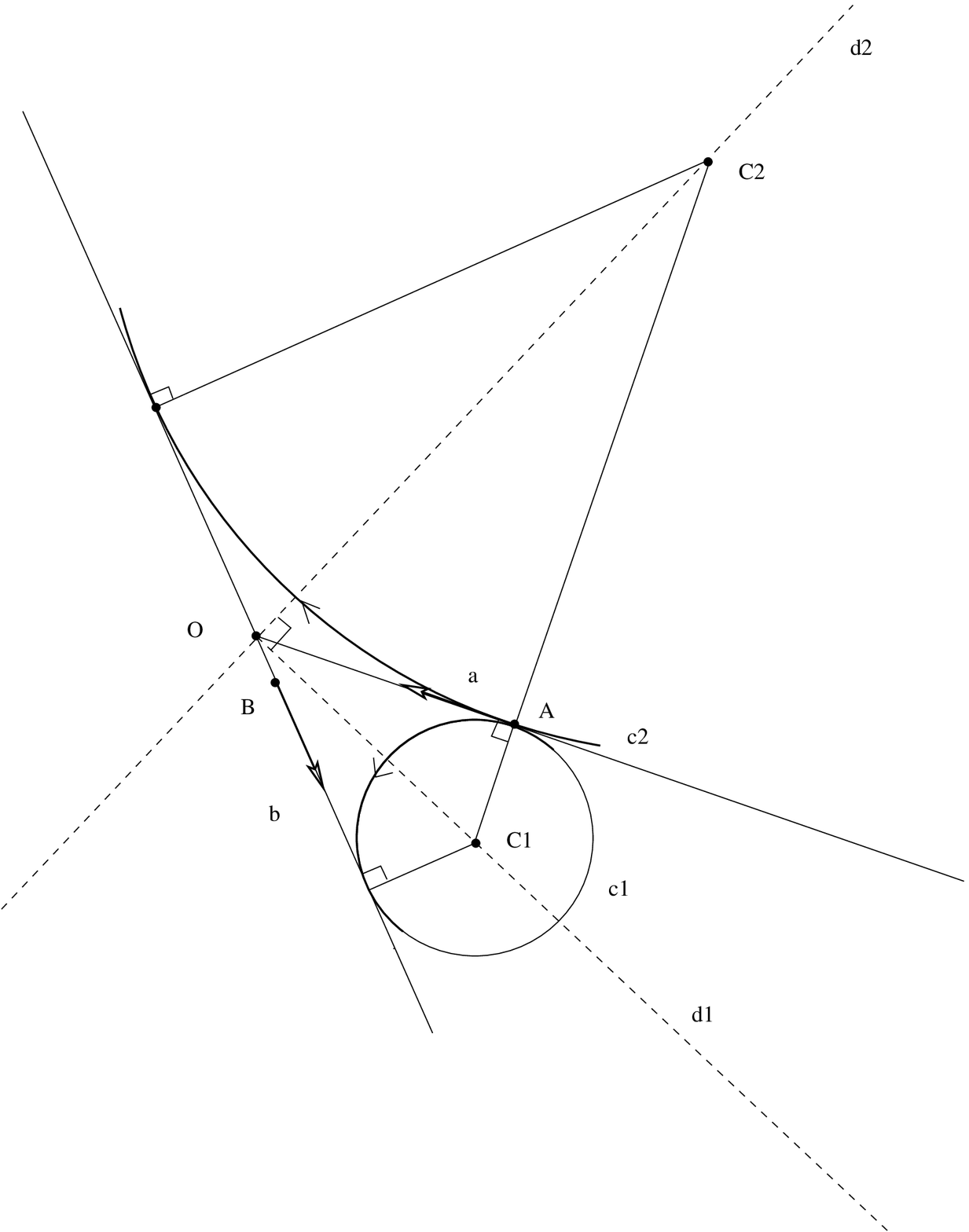, width= 10 cm} 
\end{center} 
\caption{\label{preuveuniquecercledroitelem01}$\mathcal{C}_1$ et $\mathcal{C}_2$ sont  tangents à la fois à $(OA)$ et à $(OB)$.\ifcase \cras \or /\textit{$\mathcal{C}_1$ and $\mathcal{C}_2$ are tangent to $(OA)$ and to  $(OB)$.}\fi}
\end{figure}
Supposons tout d'abord qu'il existe une courbe $Y$ de  $\mathcal{E}$
formée d'un arc de cercle $\mathcal{C}$ de longueur non nulle et d'un segment de droite $S$ éventuellement réduit à un point  et montrons qu'elle est nécessairement unique.
D'après \eqref{eq100tot}, l'arc de cercle est tangent en $A$ à $(OA)$ ou tangent en $B$ à $(OB)$. 

\begin{enumerate}

\item
Supposons  que 
$\mathcal{C}$  est tangent en $A$ à $(OA)$. Dans ce cas, le segment de droite  est  tangent en $B$ à $(OB)$ et il est donc inclus  dans la droite $(OB)$.
Puisque la courbe est de classe ${\mathcal{C}}^1$, $\mathcal{C}$ et $S$ sont tangents donc $\mathcal{C}$ est tangent à $(OB)$. 
Il n'existe que tels deux arcs de cercles possibles $\mathcal{C}_1$ et $\mathcal{C}_2$, dont les centres sont respectivement 
sur la bissectrice $d_1$ de $\widehat{AOB}$ ou la droite $d_2$ perpendiculaire à cette bissectrice, passant par $O$.
Voir figure \ref{preuveuniquecercledroitelem01}. Le centre $C_1$  (resp. $C_2$) $\mathcal{C}_1$ (resp. de  $\mathcal{C}_2$)   se trouve nécessairement sur la droite passant par $A$ et 
$d_1$ (resp. $d_2$), qui se coupent nécessairement, compte tenu de \eqref{eq20}.
Ces deux arcs de cercles, de part et d'autre de $(OA)$ doivent être parcourus dans le sens trigonométrique, pour respecter \eqref{eq140}
et, sur la figure  \ref{preuveuniquecercledroitelem01}, seul l'arc de cercle  $\mathcal{C}_1$ respecte ce sens.
Ainsi, si l'arc de cercle est tangent en $A$ à $(OA)$, il est unique et correspond à l'arc de cercle  $\mathcal{C}_1$ de la figure \ref{preuveuniquecercledroitelem01}. 
\begin{figure}[h]    
\psfrag{A}{$A$}
\psfrag{B}{$B$}
\psfrag{O}{$O$}
\psfrag{E}{$E$}
\psfrag{F}{$F$}
\psfrag{a}{$\vec \alpha$}
\psfrag{b}{$\vec \beta$}
\psfrag{Ca}{$\mathcal{C}_1$}
\psfrag{Cb}{$\mathcal{C}_3$}
\begin{center} 
\epsfig{file=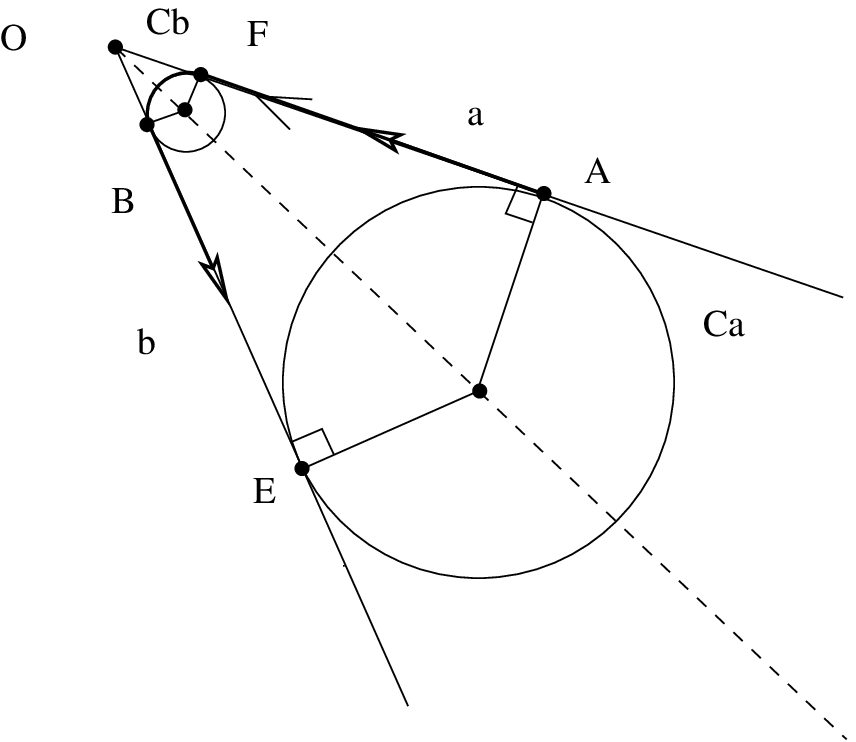, width= 6 cm} 
\end{center} 
\caption{\label{preuveuniquecercledroitelem02}Construction de $\mathcal{C}=\mathcal{C}_3$, cas 1.\ifcase \cras \or /\textit{Construction of $\mathcal{C}=\mathcal{C}_3$, case  1.}\fi}
\end{figure}
On considère ensuite le point de contact $E$ entre  $\mathcal{C}_1$, nécessairement sur la demi-droite $]OB)$. Voir figure \ref{preuveuniquecercledroitelem02}.
Si $E$ est strictement plus loin de $B$ que $O$ (ce qui implique $OB<OA$), alors la courbe $Y$ de  $\mathcal{E}$ ne peut revenir à $B$ par un segment de droite. Nécessairement,
$\mathcal{C}$  est tangent en $B$ à $(OB)$ et on construit l'unique arc de cercle $\mathcal{C}_3$, comme précédemment (voir figure \ref{preuveuniquecercledroitelem02}), tangent 
à $(OB)$ en $B$ et tangent à $(OA)$.
 Ainsi, dans ce cas, la courbe $Y$ est nécessairement formée de $\mathcal{C}=\mathcal{C}_3$ et du segment $[AF]$, non réduit à un point.
\begin{figure}[h]    
\psfrag{A}{$A$}
\psfrag{B}{$B$}
\psfrag{O}{$O$}
\psfrag{E}{$E$}
\psfrag{F}{$F$}
\psfrag{a}{$\vec \alpha$}
\psfrag{b}{$\vec \beta$}
\psfrag{Ca}{$\mathcal{C}_1$}
\begin{center} 
\epsfig{file=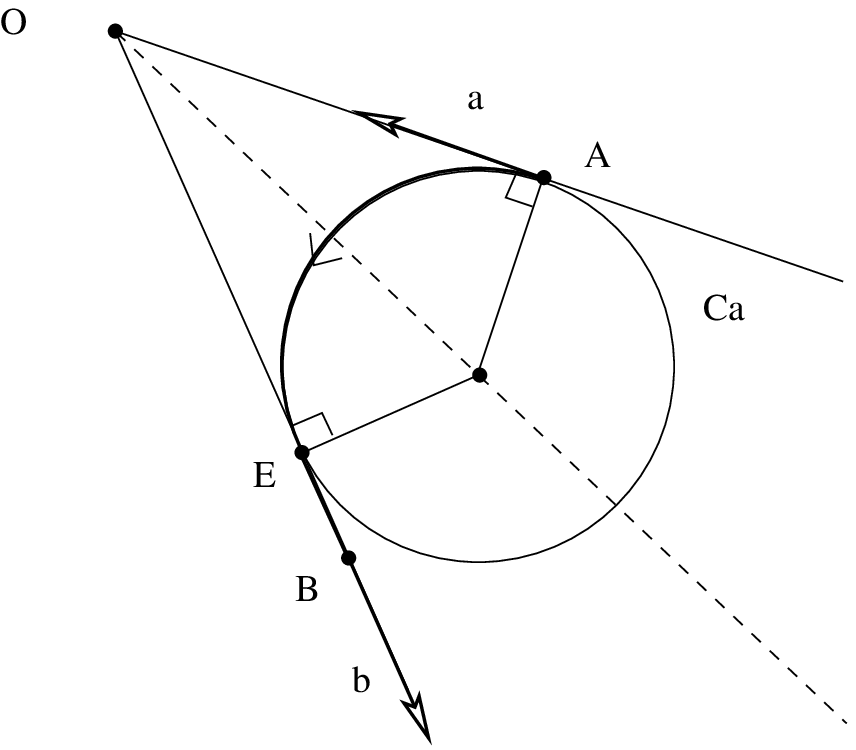, width= 6 cm} 
\end{center} 
\caption{\label{preuveuniquecercledroitelem03}Construction de $\mathcal{C}=\mathcal{C}_1$, cas 2.\ifcase \cras \or /\textit{Construction of $\mathcal{C}=\mathcal{C}_1$, case  2.}\fi}
\end{figure}
Si, au contraire $E$ est strictement moins loin de $B$ que $O$ (ce qui implique $OB>OA$), alors, de même,  la courbe $Y$ est nécessairement formée de $\mathcal{C}=\mathcal{C}_1$ et du segment $[EB]$, non réduit à un point (voir figure \ref{preuveuniquecercledroitelem03}).
\begin{figure}[h]    
\psfrag{A}{$A$}
\psfrag{B}{$B$}
\psfrag{O}{$O$}
\psfrag{a}{$\vec \alpha$}
\psfrag{b}{$\vec \beta$}
\psfrag{Ca}{$\mathcal{C}_1$}
\begin{center} 
\epsfig{file=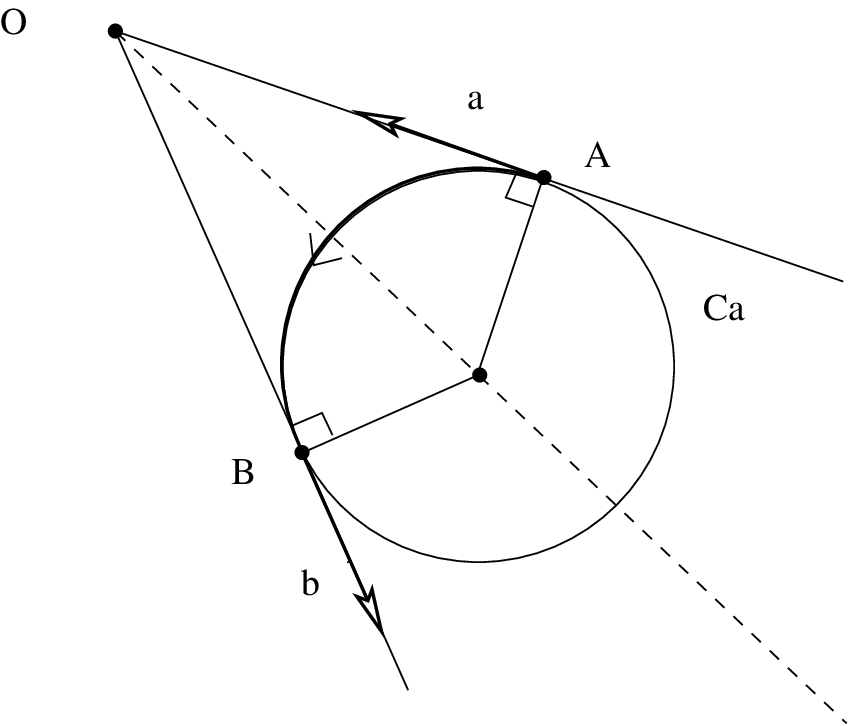, width= 6 cm} 
\end{center} 
\caption{\label{preuveuniquecercledroitelem04}Construction de $\mathcal{C}=\mathcal{C}_1$, cas 3 (cas symétrique).\ifcase \cras \or /\textit{Construction of $\mathcal{C}=\mathcal{C}_1$, case  1, symmetric  case.}\fi}
\end{figure}
Enfin, si $E$ est exactement aussi loin de $B$ que $O$ (ce qui implique $OB=OA$), alors, de même,  la courbe $Y$ est nécessairement uniquement formée de $\mathcal{C}=\mathcal{C}_1$  (voir figure \ref{preuveuniquecercledroitelem04}).
S'il existe une courbe $Y$ de  $\mathcal{E}$
formée d'un arc de cercle $\mathcal{C}$ de longueur non nulle, alors on tombe sur ce dernier cas.

\item
Si 
$\mathcal{C}$  est tangent en $B$ à $(OB)$, on arrive à la même construction. 
\end{enumerate}

Montrons maintenant l'existence de la courbe.
Dans les trois cas évoqués dans l'unicité (selon que $OB>OA$, $OB<OA$ ou $OA=OB$), on peut construire, une courbe 
$Y$ de $\mathcal{E}$,
formée d'un arc de cercle et de
périmètre dans $]0,R_a(O,A,B) \pi [$  et d'un segment de droite. 
Dans cette construction, on vérifie que 
\eqref{eq100b},
\eqref{eq100c},
\eqref{eq100d}
et 
\eqref{eq100e}
on lieu  grâce à \eqref{eq10tottot}.

L'égalité \eqref{uniquecercledroitelem01eq01} est triviale, puisque 
$\vnorm{X''}=0$ ou $1/R_a$.

\ifcase \cras
\end{proof}
\or
\end{preuve}
\fi


\ifcase \cras
\begin{proof}[Démonstration  du lemme \ref{optdubinslem01}]\
\or
\begin{preuve}[ du lemme \ref{optdubinslem01}]\
\fi

\ifcase \cras
\label{optdubinslem01preuve}
\fi

Soit $R_a=R_a(O,A,B)$, le nombre défini dans le lemme \ref{uniquecercledroitelem01}.
Supposons par exemple que $OB<OA$. On est donc dans le premier cas de la démonstration du lemme \ref{uniquecercledroitelem01}.
La courde Dubins  $\mathcal{G}(R)$ pour $R<R_a(O,A,B)$ est unique et est formée alors de deux arcs de cercles,
reliés par un segment de droite $[ED]$ avec $E\not =D$, comme le montre la figure \ref{preuve}.\ref{codu0bis}.
Le premier arc de cercle est nécessairement tangent à la droite $(OB)$ en $B$ et son centre $F$, est tel que $(BF)$ est 
$(OB)$ soient perpendiculaires. Puisque $R<R_a(O,A,B)$, $F$ est dans le secteur de plan défini par les demi-droites
$[OB)$ et $d$, la bissectrice de l'angle $\widehat{AOB}$. Si $G$ est le centre du second cercle de la courbe
de Dubins, tangent en $A$ à $(OA)$, alors $FGDE$ est un rectangle et le segment $[ED]$ est inclus 
dans le secteur de plan défini par les demi-droites
$[OB)$ et $[OA)$. Le second cercle de la courbe de Dubins est de longueur non nulle. On vérifie que la courbe $\mathcal{G}(R)$ appartient bien à $\mathcal{E}$. Cela est vrai tant que 
$R$ appartient à $]0,R_a(O,A,B)[$. Si $R$ croît, $ED$ augmente, la longueur du premier arc de cercle augmente, celle du second diminue.
Pour le cas limite, $R=R_a(O,A,B)$, correspondant  à la figure  \ref{preuve}.\ref{coducl1bis},
le premier cercle devient tangent à $(OA)$ en $E$, $D$ se confond avec $A$ et la longueur du second cercle
devient nulle. On retrouve donc l'unique courbe du lemme \ref{uniquecercledroitelem01}.
\newcommand{\lengthdiftypcodusc}{7}
\begin{figure}
\psfrag{A}{$A$}
\psfrag{B}{$B$}
\psfrag{C}{$C$}
\psfrag{D}{$D$}
\psfrag{O}{$O$}
\psfrag{D}{$D$}
\psfrag{E}{$E$}
\psfrag{F}{$F$}
\psfrag{G}{$G$}
\psfrag{a}{$\vec \alpha$}
\psfrag{b}{$\vec \beta$}
\centering
\subfigure[\label{codu0bis}Le cas 0 : $0<R<R_a(O,A,B)$\ifcase \cras \or /\textit{{The case 0 : $0<R<R_a(O,A,B)$.}}\fi]
{\epsfig{file=codu0, width=\lengthdiftypcodusc  cm}}
\quad
\subfigure[\label{coducl1bis}Le cas limite 0-1 : $R=R_a(O,A,B)$\ifcase \cras \or /\textit{{The limit case 0-1 : $R=R_a(O,A,B)$.}}\fi]
{\epsfig{file=coducl1, width=\lengthdiftypcodusc   cm}}
\quad
\subfigure[\label{codu2}Le cas  1 : $R_a(O,A,B)<R<R_b(O,A,B)$\ifcase \cras \or /\textit{{The case 1 : $R_a(O,A,B)<R<R_b(O,A,B)$.}}\fi]
{\epsfig{file=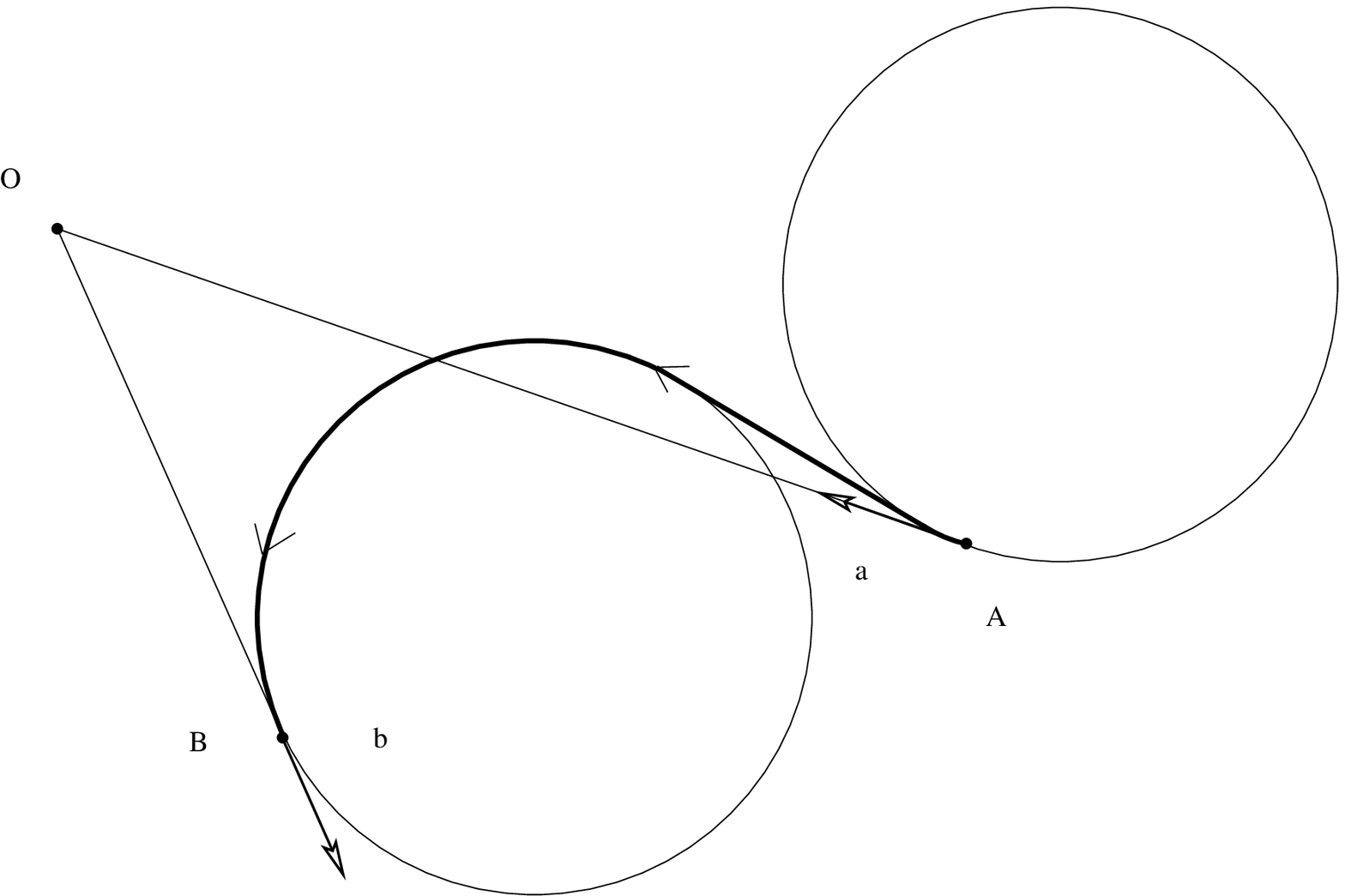, width=\lengthdiftypcodusc   cm}}
\quad
\subfigure[\label{coducl2}Le cas limite 1-2 : $R=R_b(O,A,B)$\ifcase \cras \or /\textit{{The limit case 1-2 : $R=R_b(O,A,B)$.}}\fi]
{\epsfig{file=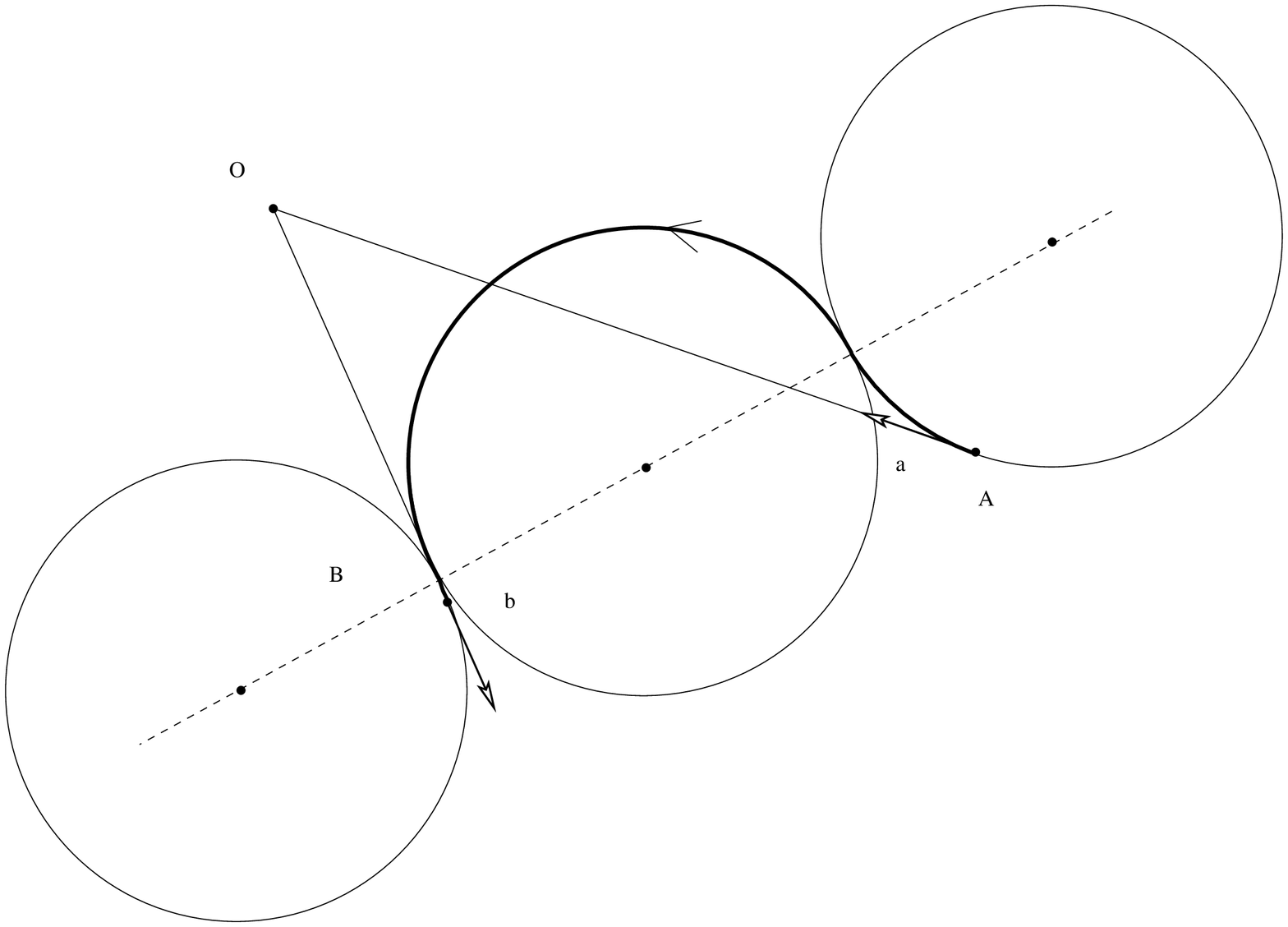, width=\lengthdiftypcodusc   cm}}
\quad
\subfigure[\label{codu3}Le cas 2 : $R>R_b(O,A,B)$\ifcase \cras \or /\textit{{The case 2 : $R>R_b(O,A,B)$.}}\fi]
{\epsfig{file=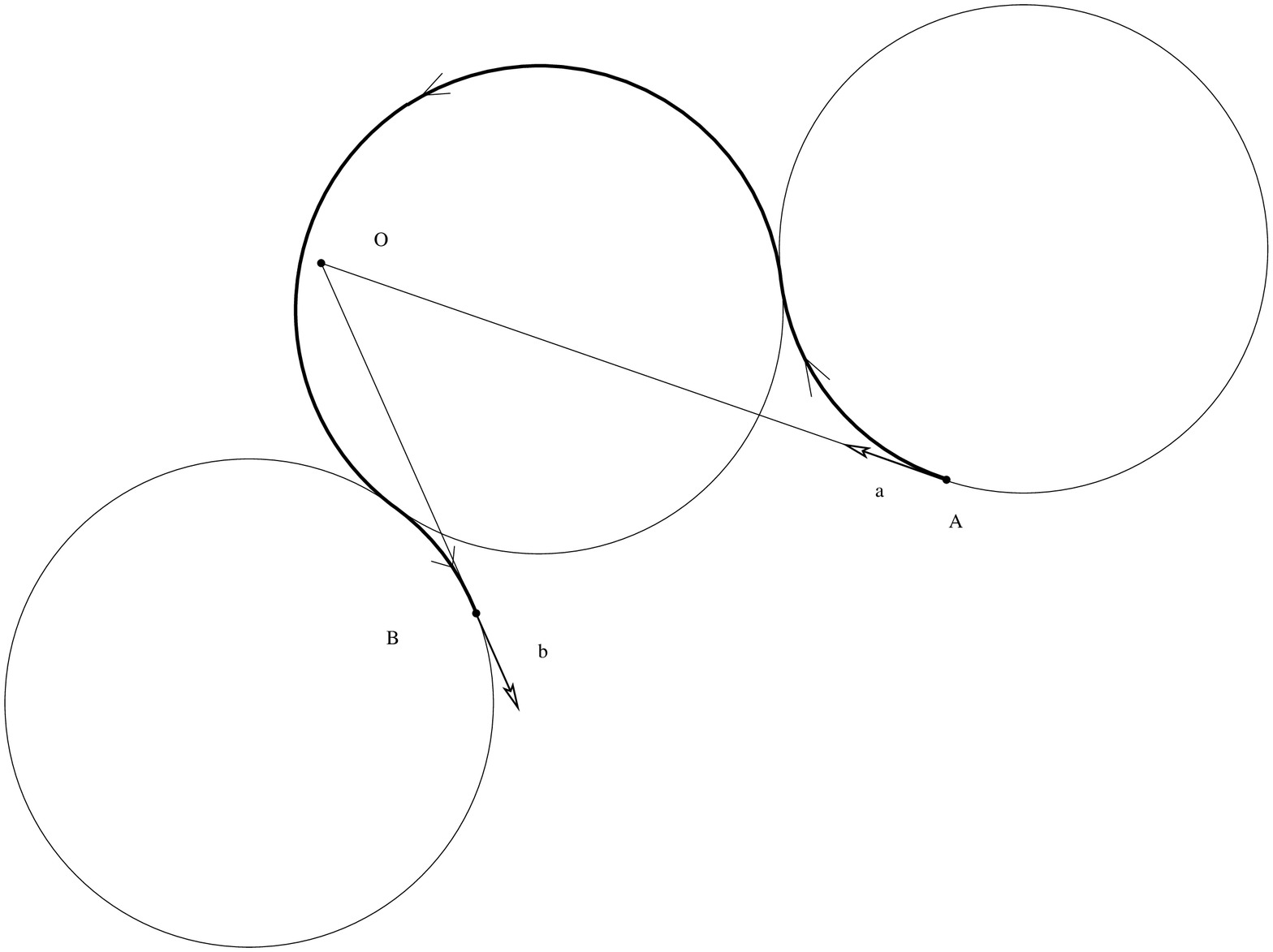, width=\lengthdiftypcodusc   cm}}
\caption{\label{diftypcodubis}Les différents types de courbes de Dubins définies par $R$  dans le cas non symétrique. $R_b(O,A,B)$ ne dépend que de $O$, $A$ et $B$. \ifcase \cras \or /\textit{{The different kinds of Dubins's Curve defined by $R$ in the non symmetric case. $R_b(O,A,B)$ depends only on $O$, $A$ and  $B$.}}\fi}
\end{figure}
\begin{figure}
\psfrag{A}{$A$}
\psfrag{B}{$B$}
\psfrag{C}{$C$}
\psfrag{D}{$D$}
\psfrag{O}{$O$}
\psfrag{a}{$\vec \alpha$}
\psfrag{b}{$\vec \beta$}
\centering
\subfigure[\label{codu0sc}Le cas 0 : $0<R<R_a(O,A,B)$.\ifcase \cras \or /\textit{{The case 0 : $0<R<R_a(O,A,B)$.}}\fi]
{\epsfig{file=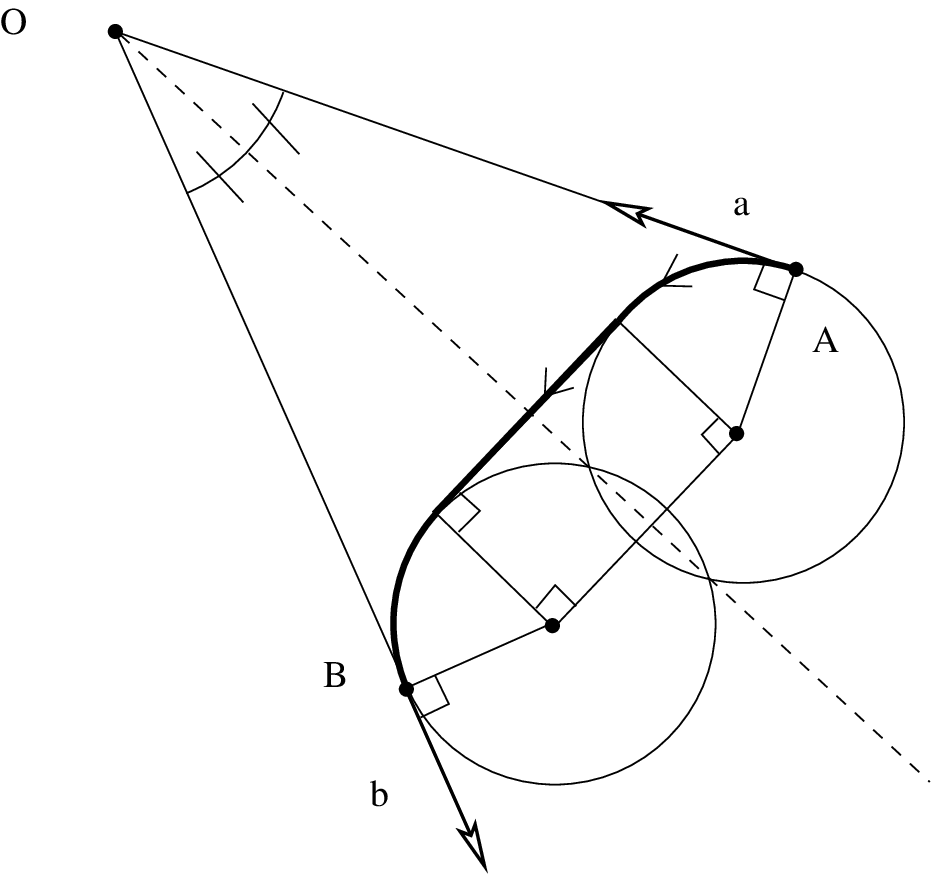, width=\lengthdiftypcodusc cm}}
\quad
\subfigure[\label{coducl1sc}Le cas limite 0-1 : $R=R_a(O,A,B)$.\ifcase \cras \or /\textit{{The limit case 0-1 : $R=R_a(O,A,B)$.}}\fi]
{\epsfig{file=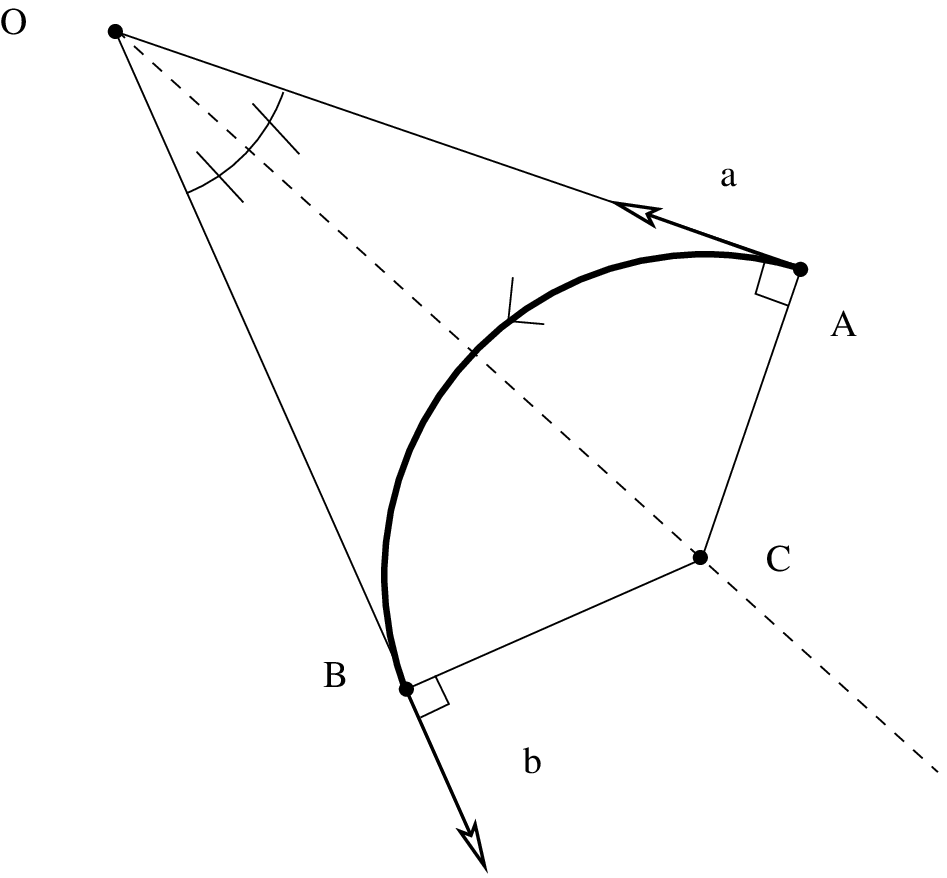, width=\lengthdiftypcodusc  cm}}
\quad
\subfigure[\label{codu2sc}Le cas  1 : $R_a(O,A,B)>R$. \ifcase \cras \or /\textit{{The case 1 : $R>R_a(O,A,B)$.}}\fi]
{\epsfig{file=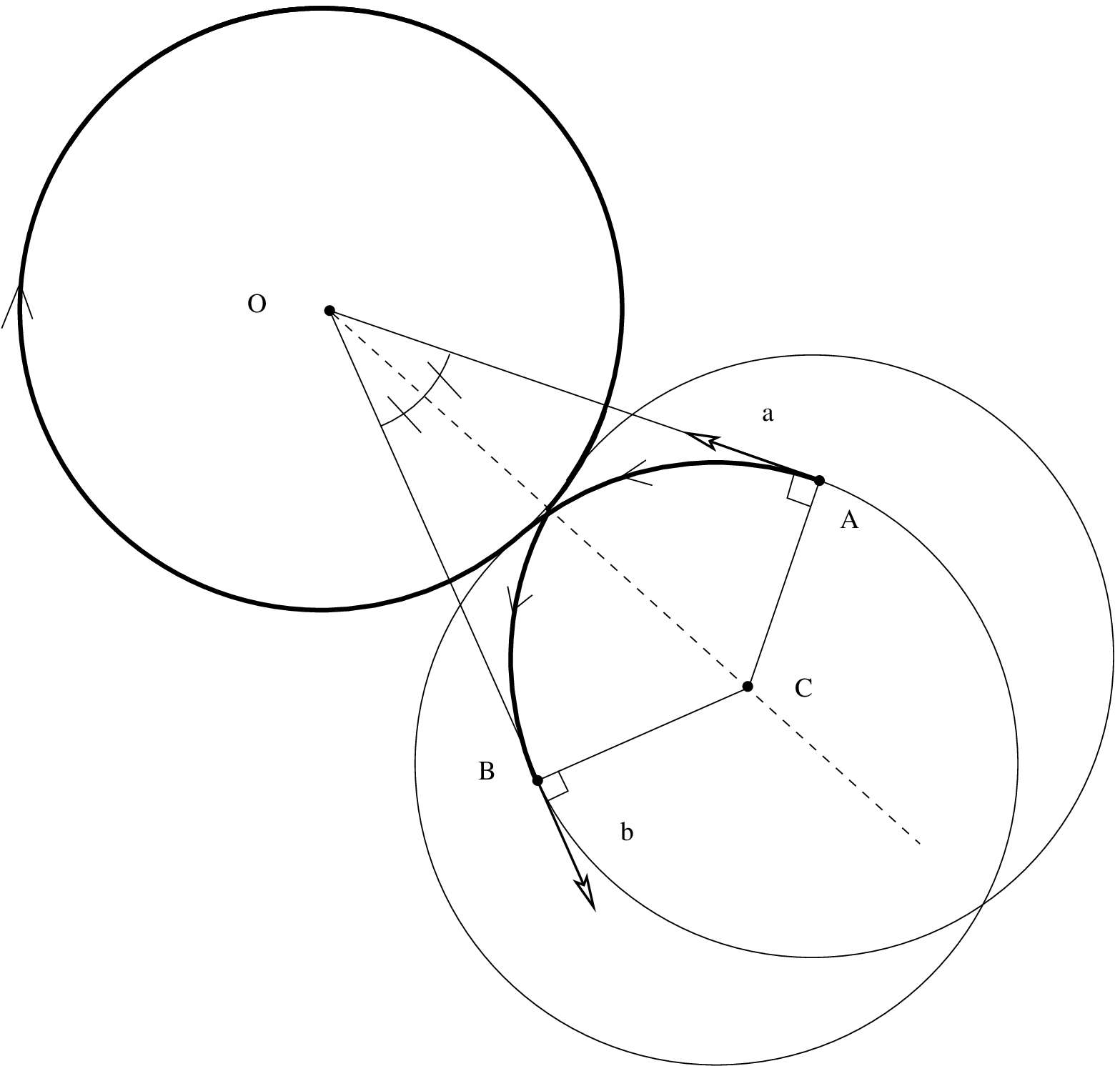, width=\lengthdiftypcodusc  cm}}
\caption{\label{diftypcodusc}Les différents types de courbes de Dubins définies par $R$  dans le cas symétrique.\ifcase \cras \or /\textit{{The different kinds of Dubins's Curve defined by $R$ in the symmetric case.}}\fi}
\end{figure}

Les cas $OA<OB$ ou $OA=OB$ se traitent de la même façon (voir figures
\ref{preuve}.\ref{codu0sc} et 
\ref{preuve}.\ref{coducl1sc}).

En revanche, on laisse au lecteur le soin de vérifier que 
  dès que $R$ dépasse strictement $R_a(O,A,B)$ 
on obtient des courbes qui présentent un changement strict du signe de la courbure. 
Voir figures \ref{diftypcodubis} et \ref{diftypcodusc}.
Le cas limite des  figures
\ref{preuve}.\ref{coducl1bis} et \ref{preuve}.\ref{coducl1sc} correspond exactement à la courbe décrite dans le lemme \ref{uniquecercledroitelem01}.
\ifcase \cras
\end{proof} 
\or
\end{preuve} 
\fi


\ifcase \cras

\ifcase \cras
\begin{proof}[Démonstration du théorème \ref{existencetheoprop01discprop01}]\
\or
\begin{preuve}[  du théorème \ref{existencetheoprop01discprop01}]\
\fi

\ifcase \cras
\label{existencetheoprop01discprop01preuve}
\fi

Il est nécessaire de prendre $p\geq 2$, car dans le cas $p=1$, le système linéaire 
 est sous-déterminé.

Tout d'abord, constatons que  ${\mathcal{E}}_p$ est inclus dans ${\mathcal{E}}$. Ainsi, 
d'après \eqref{existencetheolem02dfg}
on peut affirmer 
qu'il  existe $\delta_p>0$ telle que pour toute courbe $X$ de
 ${\mathcal{E}}_p$,
\begin{equation}
\label{discretcaspthenconeq50}
\vnorm[L^{\infty}(0,L)]{\vnorm{X''}}\geq \delta_p.
\end{equation}
On laisse vérifier  au lecteur que $X$ appartient à ${\mathcal{E}}_p$ ssi $X$ appartient à ${\left(\Erps\right)}^{p}$ et vérifie 
le système linéaire 
$\mathcal{A}X=\mathcal{B}$.
Si ${\mathcal{E}}_p$ est non vide,
il existe donc $\delta_{p,0}>0$ défini par 
\begin{equation}
\label{discretcaspthenconeq60}
\delta_{p,0}=
\inf_{X\in {\mathcal{E}}_p} \vnorm[L^{\infty}(0,L)]{\vnorm{X''}}.
\end{equation}
Concluons en 
utilisant \eqref{lemmeBeq01}.
D'après \eqref{discretcaspthenconeq60}, il existe une suite $X^n$ de ${\left(\Erps\right)}^{p}$ telle que
\begin{equation}
\label{discretcaspthenconeq70}
\lim_{n\to \infty} \vnorm[L^{\infty}(0,L(X))]{\vnorm{(X^n)''}}=\delta_{p,0}.
\end{equation}
Montrons que $X^n$ admet une suite extraite convergeant vers un élément $X$ de ${\mathcal{E}}_p$,
ce qui nous permettra de conclure.
Remarquons que  ${\mathcal{E}}_p$ est inclus dans un compact de $\Er^p$. En effet, si $X$ appartient à 
${\left(\Erps\right)}^{p}$, on a
\begin{equation*}
\forall k\in \{1,...,p\},\quad R_k>0,
\end{equation*}
et \textit{a fortiori}
\begin{equation}
\label{discretcaspthenconeq80}
\forall X=(R_1,...,R_p)\in {\mathcal{E}}_p,\quad
\forall k\in \{1,...,p\},\quad R_k\geq 0.
\end{equation}
De plus, puisque ${\mathcal{E}}_p$ est inclus dans $\mathcal{E}$, \eqref{lemmeBeq01} s'applique :
\begin{equation}
\label{discretcaspthenconeq90}
\forall X\in {\mathcal{E}}_p,\quad
L(X)\leq \nu.
\end{equation}
Dans le cas discret, on a, en notant $X=(R_1,...,X{p})$, 
\begin{equation*}
L(X)=\sum_{k=1}^{p} R_k \theta_0=
\theta_0 \sum_{k=1}^{p} R_k,
\end{equation*}
et donc, d'après \eqref{discretcaspthenconeq80} et \eqref{discretcaspthenconeq90}, 
\begin{equation}
\label{discretcaspthenconeq100}
\forall X=(R_1,...,R_p)\in {\mathcal{E}}_p,\quad
\forall k\in \{1,...,p\},\quad
R_k\leq \frac{\nu}{\theta_0}.
\end{equation}
D'après \eqref{discretcaspthenconeq80}
et \eqref{discretcaspthenconeq100}, ${\mathcal{E}}_p$ est inclus dans un fermé borné de $\Er^p$, qui est donc compact. Il existe donc $X\in\Er^p$ et  une sous-suite de $X^n$, notée de la même façon, telle que 
\begin{equation}
\label{discretcaspthenconeq110}
\lim_{n\to \infty} X^n=X.
\end{equation}
Puisque, pour tout $n$, on a 
 $X^n\in {\mathcal{E}}_p$, on a 
\begin{equation*}
\mathcal{A}X^n=
\mathcal{B},
\end{equation*}
et par continuité :
\begin{equation}
\label{discretcaspthenconeq120}
\mathcal{A}X=
\mathcal{B}.
\end{equation}
Enfin, d'après 
\eqref{discretcaspthenconprop01eq40}
 et 
\eqref{discretcaspthenconeq70}, on a, en notant $X^n=(R_1^n,...,R_p^n)$, 
\begin{equation*}
\lim_{n\to \infty}
\frac{1}{\displaystyle{\min_{1\leq k\leq p} R^n_k}}=\delta_{p,0}.
\end{equation*}
et en particulier, il existe $N$ tel que 
\begin{equation*}
\forall n\geq N,\quad 
\frac{1}{\displaystyle{\min_{1\leq k\leq p} R^n_k}}\leq 2\delta_{p,0},
\end{equation*}
ce qui implique 
\begin{equation*}
\forall n\geq N,\quad 
\min_{1\leq k\leq p} R^n_k\geq \frac{1}{2\delta_{p,0}},
\end{equation*}
et donc
\begin{equation*}
\forall n\geq N,\quad 
\forall k\in\{1,...,p\},\quad
R^n_k\geq \frac{1}{2\delta_{p,0}},
\end{equation*}
et à la limite 
\begin{equation}
\label{discretcaspthenconeq130}
\forall k\in\{1,...,p\},\quad
R_k>0.
\end{equation}
D'après \eqref{discretcaspthenconeq120} et \eqref{discretcaspthenconeq130}, on a donc $X\in{\mathcal{E}}_p$.
\ifcase \cras
\end{proof}
\or
\end{preuve} 
\fi

\fi

\ifcase \cras
\or

\clearpage

\section{Preuves complètes des résultats essentiels de cette Note.}
\label{preuvecomplete}

\input{vrai_probleme_pointI}

\fi

\end{document}